\documentclass[twoside,11pt]{article}

\title{Pseudohomology and Homology }
\author{Peter J. Kahn}

\newcommand{\Ps}[1]{\ensuremath{\Psi_{#1}}}
\newcommand{\Hm}[1]{\ensuremath{H_{#1}}}
\newcommand{\interpar}{\vspace{.2in}}
\newcommand{\sint}[1]{\ensuremath{\stackrel{\circ}{#1}}}
\newcommand{\pri}{\ensuremath{^{\prime}}}

\newcommand{\up}[1]{\mbox{\raisebox{.6ex}{#1}}}

\newcommand{\ld}[1]{\ensuremath{\ell d(#1)}}

\newcommand{\C}{\ensuremath{\mathcal{C}}}

\newtheorem{theorem}{Theorem}

\newtheorem{defn}{Definition}
\newtheorem{lemma}{Lemma}[section]
\newtheorem{cor}[lemma]{Corollary}
\newtheorem{prop}{Proposition}
\newtheorem{corr}[prop]{Corollary}

\begin{document}
\maketitle

\begin{abstract}
The notion of a {\it pseudocycle} is introduced
in ~\cite{ms} to provide a framework for defining Gromov-Witten invariants and quantum
cohomology.  This paper studies the bordism groups of pseudocycles, called
{\it pseudohomology groups}.  These satisfy the Eilenberg-Steenrod axioms, and, for smooth
compact manifold pairs, pseudohomology is naturally equivalent to homology.  Easy
examples show that this equivalence does not extend to the case of non-compact manifolds.
\end{abstract}

\section{Introduction}

Moduli spaces of pseudoholomorphic curves in a symplectic manifold $X$
are generically smooth, finite-dimensional manifolds, but usually they are not compact. 
Thus, fundamental cycles are not available for extracting homology information or for defining
invariants of $X$.  One solution to this defect is to use a theorem of Gromov~\cite{g} to
compactify the moduli spaces (e.g., see~\cite{rt},~\cite{s}).  An alternative approach, as
presented in M\up{c}Duff and Salamon~\cite{ms},  also applies Gromov's ideas by making
use of  certain dimension restrictions on the limit sets of evaluation maps associated
with the moduli spaces.    M\up{c}Duff and Salamon   call smooth
$X$-valued maps satisfying these restrictions \emph{pseudocycles} in $X$ and use them to
provide a framework for defining the Gromov-Witten invariants and quantum
cohomology of $X$.  This paper studies the bordism groups of pseudocycles in
$X$, which we call the \emph{pseudohomology} groups of $X$.  To state our first result,
let $\C$ denote the category of all smooth manifold pairs and smooth maps, as described
in \S 3,  let $\C^{prop}$ be the subcategory with the same objects and with morphisms the
\emph{proper} smooth maps, and let $\mathcal{A}$ denote the category of abelian groups
and homorphisms. Pseudohomology groups are defined for all pairs in $\C$.
\interpar

\begin{theorem}
 Pseudohomology groups and map-composition define functors
$\Ps{k}:\C^{prop}\rightarrow\mathcal{A}$, for all
$k=0,1,2,\ldots\: $.  These
functors satisfy the axioms of Eilenberg and Steenrod: homotopy-invariance,
exact-sequence, excision, and dimension.  In particular, for the dimension axiom, we
have a 
canonical isomorphism
$\Hm{k}(pt)\rightarrow\Ps{k}(pt)$, for every $k$.

\end{theorem}

\interpar

The uniqueness theorem of Eilenberg and Steenrod, which they prove over the category
of finite simplicial pairs, suggests that the isomorphism
$H_{\ast}(pt)\rightarrow\Ps{\ast}(pt)$ should extend purely formally to a natural
equivalence of functors
$H_{\ast}\rightarrow\Ps{\ast}$ over a suitable subcategory of $\mathcal{C}^{prop}$. In
fact, this is the case, but we prefer to present a concrete realization of this
equivalence.

In ~\cite{ms}, Remark  , M\up{c}Duff and Salamon outline a construction in the absolute
case that provides a candidate for such a realization.   Indeed, this paper was originally
motivated by an effort to extend that construction and understand its properties, and to
do some computations.  Our next result shows that the construction extends to the case of
pairs and does realize the equivalence that we want. In fact, we show not only that it
gives an equivalence but that the equivalence is essentially unique.

\begin{theorem}

 There
is a construction applying to all pairs
$(X,A)$ in
$\mathcal{C}$, 
\[\{\mbox{relative k-cycles in}\:\: (X,A)\}\Longrightarrow\
\{\mbox{relative k-pseudocycles in}\:\: (X,A)\},\] with the following properties:
\begin{enumerate}
  \item It induces well-defined homomorphisms
$\psi:\Hm{\ast}(X,A)\rightarrow\Ps{\ast}(X,A)$, which commute with boundary/connecting
homomorphisms.

  \item Restricted to $\mathcal{C}^{prop}$, the maps $\psi$ form a natural
transformation of functors $\Hm{\ast}\rightarrow\Ps{\ast}$.

  \item  $\psi:\Hm{\ast}(pt)\rightarrow\Ps{\ast}(pt)$ equals the isomorphism of
Theorem 1. 

  \item Let $\mathcal{C}^{prop}_{hand}$ denote the full subcategory of
$\mathcal{C}^{prop}$ with objects the pairs $(X,A)$ admitting finite handle decompositions
(cf. \S 3).  Then the natural transformation $\psi$ is an \emph{equivalence} of functors
on
$\mathcal{C}^{prop}_{hand}$.

Furthermore, $\psi$ is the only natural transformation $H_{\ast}\rightarrow \Psi_{\ast}$
on $\mathcal{C}^{prop}_{hand}$ satisfying a and c.  
 
\end{enumerate}

\end{theorem}

\interpar

\noindent{\bf Remarks:}\begin{enumerate}

  \item There are simple examples (Examples 1,2
\S 4) which show that although 
$\Ps{\ast}$ is defined for the objects of $\mathcal{C}$, it does not extend functorially
to the maps in $\mathcal{C}$.  
 
  \item The notion of a finite handle decomposition of a pair $(X,A)$ as described in \S3
is a relative notion.  Essentially, $X$ is obtained from $A$ by attaching finitely
many handles.  Thus, if $X=(X,\emptyset)$ is a non-compact manifold in
$\mathcal{C}^{prop}$, it is not an object of $\mathcal{C}^{prop}_{hand}$.  

Theorem 2 d does not
extend to
$\mathcal{C}^{prop}$, which is to say that there is a substantial deviation between
homology and pseudohomology for non-compact manifolds, as the following two results
indicate:
		\begin{enumerate}
    \item \begin{prop}
								Let $X$ be a non-compact, connected, orientable manifold with empty boundary.  There
is no surjective homomorphism $\Hm{\ast}(X)\rightarrow\Ps{\ast}(X)$.
          \end{prop}
    \item \begin{prop}
								Let $X$ be a smooth, compact, orientable manifold.  There is no natural transformation
$H_{\ast}\rightarrow\Ps{\ast}$ on $\mathcal{C}^{prop}$ that maps $H_{\ast}(X\times
R)\rightarrow\Ps{\ast}(X\times R)$ injectively.
           \end{prop}
    \end{enumerate}
  \item  Unlike homology, $\Ps{\ast}$ does not commute with (infinite) direct limits.  For
if it did, neither of the above propositions would hold.  Of course, this can be seen
even more easily in the example of the identity map of a non-compact $0$-manifold $X$. 
Such a map is a $\Psi_{0}$-cycle which clearly is not in the subgroup of $\Ps{0}(X)$
generated by the pseudohomology of the compact submanifolds of $X$.  

  \item  It is worth pointing out the connection between classical bordism and
pseudocycle bordism:

\emph{ Let
$\Omega_{\ast}(X,A)$ denote the
classical, oriented bordism group of
$(X,A)$, which we may take to be defined via smooth maps,  and let
$\mu:\Omega_{\ast}(X,A)\rightarrow
H_{\ast}(X,A)$ be the standard
fundamental-class-evaluation
homomorphism.  Then, $\psi\circ\mu:\Omega_{\ast}(X,A)\rightarrow\Ps{\ast}(X,A)$ is
induced by the inclusion of bordism classes.}

This result is an immediate consequence of the construction and results of \S\S 5,6.

  \item In ~\cite{ms}, M\up{c}Duff and Salamon define an intersection pairing for pseudocycles
in compact manifolds, which is analogous to the standard intersection pairing for homology. It is not hard to see that $\psi$ respects intersection pairings: i.e.,
$\alpha\cdot\beta=\psi(\alpha)\cdot\psi(\beta)$, for all homology classes $\alpha,\;\beta$ of complementary dimensions.  One way to see this is to use the relationship with classical bordism displayed in the preceding remark.  Both 
$\mu$ and the inclusion-induced map $\Omega_{\ast}\rightarrow\Psi_{\ast}$ respect intersection pairing: for the former, this is a well-known fact (cf.,~\cite{do}); for the latter, it follows from the definitions.  Now use the classical fact that the index of the subgroup $\mu(\Omega_{n}(X))$ in $H_{n}(X)$ is finite for all $n$ and all compact $X$.

Given a closed, connected, oriented $m$-manifold $X$, the pairing on pseudocycles can be used to associate with every $k$-dimensional pseudocycle $\phi$ in $X$ a homomorphism $\Ps{m-k}(X)\rightarrow Z$, depending only on the bordism class of $\phi$, which we continue to denote by the same letter. Composing this homomorphism with $\psi$, we obtain a homomorphism $H_{m-k}(X)\rightarrow Z$, which we call $\psi^{\ast}(\phi)$.  M\up{c}Duff and Salamon then say that a homology class $\alpha\in H_{k}(X)$ is \emph{weakly represented} by $\phi$ if the intersection product with $\alpha$ determines the same homomorphism as $\psi^{\ast}(\phi)$. We prefer to say that such $\alpha$ and $\phi$ are \emph{weakly related}.  This notion is used in ~\cite{ms} both because the precise relationship between homology and pseudohomology was not clear and because the stated homomorphism was the main object of interest.

In light of these comments and Theorem 2, one can see that the notion `weakly related' has the following interpretation: Namely, the isomorphism $\psi$ induces an isomorphism $\overline{\psi}:H_{\ast}(X)/tors\rightarrow \Psi_{\ast}(X)/tors$.
The relation of weak relatedness is precisely the pullback of $\overline{\psi}$ to a relation $H_{\ast}(X)\leftrightarrow \Psi_{\ast}(X)$. 
\end{enumerate}

\interpar

What is required for Theorem 2 is a
well-structured, geometric notion of ``$k$-cycle,'' as well as a relative version. Rourke and
Sanderson present such a notion in ~\cite{rs}, which they call a ``cycle.''  We shall
use the older term ``circuit''~\cite{l} for this, to avoid too many uses of the word
``cycle''.  A theorem of Rourke and Sanderson ~\cite{rs} asserts that bordism of singular
circuits is naturally equivalent to singular homology (cf., Theorem in \S 5).  Circuits
are piecewise linear (PL) objects.  Our procedure in
\S 5 will be to show that singular circuits (both absolute and relative) can be
smoothed (away from a codimension-two set) in an essentially unique way to produce
pseudocycles (both absolute and relative, respectively), and similarly for
bordisms of these.

Circuits have also been called pseudomanifolds (e.g.,~\cite{st}). We point out that
Ruan and Tian~\cite{rt} also use the term ``pseudomanifold'' in the context of
pseudoholomorphic curves:  in particular, in their analysis of moduli spaces (i.e., the
domains of the evaluation maps), which we have already alluded to.  Their notion is not
quite the same as that of Rourke and Sanderson (or the classical notion), since their
definition combines the PL and smooth categories, effectively subsuming some smoothing
results in their definition. Moreover, their codimension-two singular sets are always
assumed to be skeleta of triangulations, which is too restrictive for our purposes in
Theorem 2 (e.g., see
\S5).

\interpar

We conclude this introduction by describing the organization of the paper.  

Section 2
presents facts about limit sets of continuous maps and the limit dimension of
these.  These facts form the basis for deriving most of the properties of pseudocycles that
appear in Theorem 1. 

Section 3 describes the various categories of smooth manifolds that we use.  In
particular, it briefly describes the notion of manifold with corners, which turns out to
require some special attention in our arguments.  The Appendix discusses corners in
more detail. 

Section 4 presents proofs of the assertions in Theorems 1 and
Propositions 1 and 2.  Proposition 1 amounts to constructing a top-dimensional pseudocycle
that cannot be realized by a homology class. An amusing variant of Hirsch's proof of the
Brouwer Fixed Point Theorem accomplishes this. The most serious proof in this section is the
proof of the excision property for pseudohomology.  A key ingredient for this allows us to
construct lower dimensional pseudocycles via transversality (cf., Corollary 7, \S 2
)\footnote{As this paper was being completed, the author became aware of the revision of
~\cite{ms}, which is in preparation, and he thanks M\up{c}Duff and Salamon for making a
copy of some relevant chapters available.  A result closely similar to Corollary 7
appears in the revision.}

Section
5 begins with Rourke and Sanderson's definition of singular circuits and their theorem
that these may be used to define singular homology.  We show that singular
circuits can be replaced by so-called ``weakly piecewise smooth'' circuits. This step is
important, because singular circuits do not allow us to control limit dimension, whereas
weakly piecewise smooth circuits do.  Next we find certain codimension-two ``singular''
subpolyhedra in circuits  whose complements are PL manifolds, and we show that these
manifolds are smoothable in an essentially unique way. Similarly for relative circuits
and bordisms. The main work here involves these last, which are essential for obtaining
a transformation that is well defined, natural and commutes with the relevant connecting
homomorphisms. Finally, we show how to smooth the weakly piecewise smooth maps away from
the singular sets, so that away from the singular sets they are pseudocycles (resp.,
bordisms of pseudocycles).  This constitutes the desired transition from cycles to
pseudocycles.  

In Section 6, we piece together these steps to define the transformation
$\psi$ and conclude the proof of Theorem 2. 

The Appendix discusses manifolds and submanifolds with corners to the extent needed in
this paper.  Corners arise naturally in some of our constructions, and the usual
technique of ``straightening'' or ``rounding'' them requires us to modify maps near the
corners.  This could potentially cause us to lose control of the limit dimension of these
maps.  Thus, the main goal of the appendix is to show how, after straightening corners,
we may smooth maps while preserving their limit dimension.

\section{Limit sets and limit dimension}

All the topological spaces that we consider in this paper will be Hausdorff, second
countable, and locally-compact. We remind the reader that a continuous map
$f:X\rightarrow Y$ is
\emph{proper} if $f^{-1}(C)$ is compact for all compact $C$ in $Y$. In particular, this
holds if $X$ is compact or $f$ is a homeomorphism.  A homotopy will be called proper if
it is proper as a map.

\begin{defn}

	Let $f:X\rightarrow Y$ be a continuous map.  The \emph{limit set} $L(f)$ of
$f$ consists of all $y\in Y$ for which there is a sequence
$\{x_{n}\}$ in $X$ with no limit points such that $f(x_{n}) \rightarrow
y$. 
\end{defn}

The following proposition and corollaries summarize the basic facts that we
need about limit sets.

\begin{prop}

Let $f:X\rightarrow Y$, $g:Y\rightarrow Z$, and $f^{\prime}:X^{\prime}\rightarrow
Y^{\prime}$ be continuous. Then:
\begin{enumerate}

\item If $C \subseteq Y$ is compact and $C \cap\ L(f) = \emptyset$,
then $f^{-1}(C)$ is compact.

\item	$f$ is proper iff $L(f) = \emptyset$. 

\item If $X_{1}$ is closed in $X$, then $L(f|X_{1})\subseteq L(f)$.

\item If $X = X_{1}\cup X_{2}$, then $L(f) \subseteq 
L(f|X_{1})\cup L(f|X_{2})$, with equality if $X_{1}$ and $X_{2}$ are
closed.

\item $L(f\times f^{\prime}) = L(f)\times
\overline{f^{\prime}(X^{\prime})} \cup\ \overline{f(X)}\times
L(f^{\prime})$. 

\item $g(L(f)) \subseteq\ L(gf) \subseteq\ g(L(f))\cup L(g).$

Thus, if $g$ is proper, $L(gf)=gL(f)$.

\item If $f$ is surjective, then $L(g)\subseteq L(gf)$.  Thus, if $f$ is
surjective and proper, $L(g)=L(gf)$.  As a special case, if $K$ is compact and
$pr$ is the projection $X\times K
\rightarrow X$, then $L(f\circ pr) = L(f)$.

\item If $X$ is dense in a compact space $W$, and $f$ extends to a
continuous map $g:W\rightarrow Y$, then $L(f) = g(W\setminus X)$.

\item If $A$ is a closed subset of $Y$, then $L(f\mid f^{-1}(A))
\subseteq L(f) \cap A$.

\end{enumerate}

\end{prop}

These properties are straightforward consequences of the definition.  We leave the
proofs to the reader.

We shall call a homotopy $\Phi:Y\times[0,1]\rightarrow Z$ an \emph{isotopy} if the
associated map map $\Phi^{\prime}$
defined by $(y,t)\mapsto (\Phi(y,t),t)$ is a homeomorphism $Y\times[0,1]\rightarrow
Z\times[0,1]$. 

\begin{corr}
Every isotopy is proper.
\end{corr}
{\bf Proof:}  Let $pr:Y\times [0,1]\rightarrow Y$ be the projection map, so that
$\Phi= pr\circ\Phi^{\prime} $. Then, by Proposition 3 b,f,together with
the facts that $\Phi^{\prime}$ and $pr$ are proper,
\mbox{$L(\Phi)\subseteq pr(L(\Phi^{\prime}))\cup L(pr) = \emptyset$}.  Now apply
Proposition 3b.  This completes the proof.

\begin{corr}

Let $f$ be as in Proposition 3, and let $\Phi:Y\times[0,1]\rightarrow Z$ be a proper
homotopy.  Then \[ L(\Phi\circ (f\times id_{[0,1]})) = \Phi(L(f)\times [0,1]).\]

\end{corr}
{\bf Proof:}  Set $g=f\times id_{[0,1]}$.  Then, by Proposition 3 b,g,
\[\Phi(L(g))\subseteq L(\Phi\circ g)\subseteq \Phi(L(g))\cup L(\Phi)=\Phi(L(g)).\]
It remains to observe that $L(g) = L(f)\times [0,1]$ (Prop. 3b,e). 
This completes the proof.

\interpar

\begin{defn}
Let $f,g:X\rightarrow Y$ be two maps, and assume Y is endowed with a metric d. 
We say that f and g are \emph{equal at infinity}, written $f=_{\infty}g$, if,
for every sequence $\{x_{n}\}$in X without limit points, $\lim
d(f(x_{n}),g(x_{n}))=0$.
\end{defn}

For example, the functions $f,g:[1,\infty)\rightarrow[1,\infty)$ given by
$f(x)=x$ and $g(x)=x+1/x$ are equal at infinity when $[1,\infty)$ is given the
standard metric.  So are $f(x)=x$ and $g(x)=x+1$ when $[1,\infty)$ is given the
metric $d(x,y)=|x-y|/\max(x,y)$, but not when it has the standard metric. When
necessary, we shall specify the metric on the space $Y$. The following fact is an
immediate consequence of the definition.

\begin{lemma}
If $f=_{\infty}g$, then $L(f)=L(g)$.
\end{lemma}

In the next section, we shall discuss the various categories of smooth
manifolds and smooth maps that we use in this paper, including some facts about
manifolds with corners.  The following definition is invariant under any any
(reasonable) notion of smooth manifold and smooth map. 

\begin{defn}
	Let $Y$ be a smooth manifold, and let
$S$ be a non-empty subset of
$Y$.  The
\emph{smooth dimension} of $S$, denoted $sd(S)$, is the smallest integer $m$ for
which
$S$ is contained in the image of a smooth map $X\rightarrow Y$ with
$dimX\leq m$.  We define $sd(\emptyset)=-1$, in accord with the usual convention.

If $f:X\rightarrow Y$ is a smooth
map, then we define its
\emph{limit dimension}, written \ld{f}, by $\ld{f}=sd(L(f))$.
\end{defn}

The following corollary to Proposition 3 follows immediately from this definition:
\begin{corr}

Let $f:X\rightarrow Y$ and $f^{\prime}:X^{\prime}\rightarrow
Y^{\prime}$ be smooth, with $f\pri$ proper. Assume that $X\pri$ has dimension $n$.
Then:
\begin{enumerate}
  \item If $X_{1}$ is closed in $X$, then $\ld{f|X_{1}}\leq\ld{f}$.
  \item If $X = X_{1}\cup X_{2}$, then $\ld{f}\leq max\{\ld{f|X_{1}},\ld{f|X_{2}}\}$,
with equality when both $X_{1}$ and $X_{2}$ are closed.
  \item $\ld{f\times f\pri}\leq\ld{f}+n$.
  \item If $K$ is compact and $pr$ is the projection $X\times K\rightarrow X$, then
$\ld{f\circ pr}=\ld{f}$.
		\item If $g:Y \rightarrow Z$ is smooth, then $\ld{f}\leq \ld{gf}\leq
max(\ld{f},\ld{g})$. 
		\item If $X$ is dense in a compact, smooth manifold $W$, and $f$ extends to a
smooth map $g:W\rightarrow Y$, then $\ld{f}\leq sd(W\setminus X)$.

\end{enumerate}
\end{corr}

The final corollary to Proposition 3 plays an important role in the proof of the
excision property.  Suppose that $U$ and $V$ are  smooth manifolds, which may
have  non-empty boundaries $bU$ and $bV$  but  no corners, and suppose that
$(X,A)$ is a pair of smooth manifolds, no restrictions on boundary and corners.

\begin{corr}

Let $f:V\rightarrow X$ and
$g:U\rightarrow X$ be smooth maps, both transversal to the submanifold
$A\subseteq X$. We also assume that $f|bV$ and $g|bU$ are transversal to
$A$. Finally, suppose that
$L(f)\subseteq g(U)$. Then
$\ld{f|f^{-1}(A)}\leq dimU-codimA.$

\end{corr}
{\bf Proof:}  If $A$ has codimension zero, then the transversality hypotheses are vacuous,
and the conclusion reduces to the inequality $\ld{f|f^{-1}(A)}\leq dimU$. In this
case, the result follows immediately from Corollary 6 a and the fact that
$f^{-1}(A)$ is closed in
$V$. Now suppose that $codimA=n\geq1$.  Using Proposition 3 i, we calculate as follows:
$$L(f|f^{-1}(A))\subseteq L(f)\cap A\subseteq g(U)\cap A=g(g^{-1}(A)).$$
By transversality, $g^{-1}(A)$ is a smooth submanifold of $U$ of codimension equal to
$n$, which yields the desired result.

\section{Smooth manifolds}

All the smooth manifolds in this paper are assumed to be
finite-dimensional and second-countable. Beyond that, we need to impose
certain restrictions depending on the context.

\subsection{Various categories of smooth manifolds}
 
The constructions and theory considered in this paper require that we
distinguish among certain categories of smooth manifolds.  

To begin with, pseudocycles are, by definition, not
generally defined on compact manifolds.  Indeed, the theory requires this. However, as
described in the introduction, the theory also requires some sort of compactness or
properness condition on the target manifolds (or pairs). So, we need to distinguish the
class of manifolds allowed for the domains of pseudocycles from the class of manifolds
(or maps) allowed for the targets.

Further, in this paper, corners arise occasionally both in the targets of pseudocycles
and in some constructions involving the domains.  Moreover,
the standard procedures for dealing with corners can potentially affect
the limit dimensions of the maps that we are studying.  So
corners cannot be disregarded completely.   The reader is referred to the
Appendix where we give some basic definitions and facts about corners. 
Here, we give only a brief description.

Smooth manifolds with corners are modeled locally on open subsets of the
non-negative orthants of Euclidean spaces, by precise analogy with the
usual definitions of smooth manifolds and smooth manifolds with boundary. 
A  smooth submanifold of a  smooth manifold with corners is a
subset that inherits the structure of a smooth manifold
with corners in the usual way via the restriction of charts.  By a
\emph{smooth manifold pair}, we shall mean a smooth
manifold with corners, together with a smooth submanifold that is
\emph{closed} as a subset of $X$. However, we shall occasionally restrict
this generality.  Henceforth, when we use the term ``smooth manifold,'' we shall allow
the possibility of corners unless we rule this out explicitly.

\begin{defn}
Let $\mathcal{C}$ be the category with objects all smooth
manifold pairs. and with morphisms all  smooth maps of these
pairs. 
Let
$\mathcal{C}^{c}$ denote the full subcategory of $\mathcal{C}$ with
objects the compact pairs and
$\mathcal{C}^{prop}$ the subcategory with no restriction on objects but
with all maps required to be proper.

The category $\mathcal{S}$ is the full subcategory of $\mathcal{C}$ whose 
objects $(X,A)$ consist of manifolds $X$ and $A$ with no corners. We observe that if
$A$ is connected, then either $A$ is contained in the boundary $bX$ or  $bX\cap A=C$ is a
union of components of $bA$ with
$A$ transverse to $bX$ at $C$.   $\mathcal{S}^{c}$ will denote the full subcategory of
$\mathcal{S}$ with objects  the compact manifold pairs, and
$\mathcal{S}^{prop}$ is the subcategory with the same objects as
$\mathcal{S}$ but with all smooth maps required to be proper.
\end{defn}

\subsection{Handles and finite handle decompositions}

Our proof of Theorem 2 parallels the uniqueness
argument of Eilenberg and Steenrod~\cite{es}, which takes place in the category of
finite simplicial pairs.  The smooth manifold structures that we take to
parallel simplicial structures are finite handle decompositions.  We alert the
reader that our definition of these is slightly non-standard.  As noted
in the introduction, it is also possible by these methods to prove an
existence result for the natural transformation $\psi$.  However, this
seems to require a more elaborate notion of handle decomposition, which
we prefer to avoid in this paper.
\interpar

Let $D^{n}$ denote the unit $n$-ball and $S^{n-1}$ its boundary.  The
corresponding ball and sphere of radius $r$ will be denoted $rD^{n}$ and
$rS^{n-1}$, respectively.  Let $Y$ be a smooth $m$-manifold and
$f:S^{i-1}\times 2D^{m-i}\rightarrow Y$ a smooth embedding with image a
smooth submanifold of $Y$ contained in $bY$.  It follows easily from
the definition of `submanifold' that $f(S^{i-1}\times D^{m-i})=B$ is contained in an open face of $Y$---see the Appendix for
a definition of open face, an open submanifold of $bY$ not meeting the corners of $Y$. 
Introduce a corner at
$bB$ and glue
$Y$ to
$H=D^{i}\times D^{m-i}$ by identifying  $S^{i-1}\times D^{m-i}$ with $B$ via $f$. Again
as indicated in the Appendix, the resulting manifold can be given a differentiable
structure compatible with those of $Y$ (with corner added) and $H$.  We call this
smooth manifold \emph{the result of of attaching the
$i$-handle $H$ to $Y$} and denote it by $Y+H$. We use similar notation for
any finite set of handles disjointly attached to $Y$.

Alternatively, if $Y$ and $H$, as above, are smooth, codimension-zero
submanifolds of $Z$, such that $Y\cap H=B$ and $H\subseteq intZ$, then
we may `round the corner' $bB$ in $Z$, as indicated in the Appendix,
obtaining a smooth, codimension-zero submanifold $r(Y\cup H)$ which is
essentially the same as $Y+H$.  We use the two constructions
interchangeably.

\begin{defn}
Let $(X,A)$ be a smooth manifold pair in $\mathcal{C}^{prop}$, where $dimX=m$.  By a
\emph{finite handle decomposition} of $(X,A)$ we mean a nested sequence of
codimension-zero submanifolds
\[A\subset X_{-1}\subseteq X_{0}\subseteq\ldots\subseteq X_{m-1}\subseteq X_{m}=X\]
satisfying
\begin{enumerate}
  \item $X_{-1}$ is a tubular neighborhood of $A$.
  \item Each $X_{i}\,, i\geq 0$, admits a smooth, proper deformation retraction to
\[ X_{i-1}+H_{1}+H_{2}+\ldots ,\]
where $H_{1}, H_{2},\ldots$ is a finite collection of $i$-handles.
\end{enumerate}
\end{defn}
\interpar
We let $\mathcal{C}^{prop}_{hand}$ denote the full subcategory of $\mathcal{C}^{prop}$
whose objects are the pairs that admit finite handle decompositions.  Similarly for
$\mathcal{C}^{c}_{hand}$. 

We note that the objects of $\mathcal{C}^{c}_{hand}$ include all the objects of
$\mathcal{S}^{c}$.  Indeed, if $(X,A)$ is such a manifold pair with $A=\emptyset$, then
this is a classical result that can be proved using Morse functions (as in, say,
~\cite{m}). These techniques can be extended to the case of pairs. 

\section{Pseudohomology}

In this section we extend to the relative case the definition of
pseudocycle and bordism of pseudocycles presented in ~\cite{ms} ,
and we derive the basic properties of the resulting bordism (i.e.,
\emph{pseudohomology}) groups.  In effect, these results constitute a proof of Theorem 1
of the introduction.

\begin{defn}

Let $(X,A)$ be a smooth manifold pair in $\mathcal{C}$ and $k$ a
non-negative integer. A
$k$-dimensional \emph{relative pseudocycle} in $(X,A)$ (or
\emph{relative}\ \Ps{k}-\emph{cycle} in $(X,A)$) is a pair $(V,f)$ such
that $V$ is an oriented manifold in $\mathcal{S}$, $f$ is smooth, and the
following dimension conditions are satisfied: 
\begin{enumerate}
  \item $dimV = k$; 
  \item $\ld{f}\leq \max(-1,k-2)$;
  \item $\ld{f|bV}\leq\max(-1,k-3)$.
\end{enumerate}

\end{defn}

We use the following standard conventions. Namely, when $bV =
\emptyset$, we call $(V,f)$ an \emph{absolute} pseudocycle (or absolute
\Ps{k}-cycle) in
$(X,A)$.  And when $A$ is empty, we write $X$ instead of $(X,\emptyset)$
and omit the word ``absolute.''

If $(V,f)$ is a relative \Ps{k}-cycle in $(X,A)$, then clearly
$(bV, f|bV)$ is a \Ps{k-1}-cycle in $A$.

We now remind the reader of the usual notion of $\emph{oriented, relative
smooth bordism}$, but without the usual compactness condition (cf.,e.g., ~\cite{cf}).
Let $f_{i}:(V_{i},\partial V_{i})\rightarrow (X,A), i=0,1$, be smooth
maps of (oriented) $k$-manifolds. By a relative smooth bordism in
$(X,A)$ between
$f_{0}$ and $f_{1}$ we mean a pair $(W,F)$, such that: 
\begin{itemize}
  \item $W$ is an oriented $(k+1)$-manifold in $\mathcal{S}$;
  \item $bW$ contains the disjoint union $V_{0}\sqcup V_{1}$ as a codimension-$0$
submanifold, closed as a subset of $bW$;
  \item $F:W\rightarrow X$ is a smooth map extending
$f_{0}\sqcup f_{1}$;
  \item $F(\overline{bW\setminus(V_{0}\sqcup V_{1})})\subseteq A$.
\end{itemize}

\begin{defn}
A relative $\Ps{k}$-bordism between relative $\Ps{k}$-cycles
$(V_{i},f_{i})$ in $(X,A), i=0,1$ is a relative smooth bordism $(W,F)$
in $(X,A)$ between $(V_{0},f_{0})$ and $(V_{1},f_{1})$
as above, such that: 
\begin{enumerate}
  \item $\ld{F}\leq\max(-1,k-1)$;
  \item $\ld{F|\overline{bW\setminus(V_{0}\sqcup
V_{1})}}\leq\max(-1,k-2)$.
\end{enumerate}

\end{defn}
When the pseudocycles $(V_{i},f_{i})$ are absolute and $V_{0}\sqcup
V_{1}=bW$ , we call $(W,F)$ an
\emph{absolute bordism}.  When $A=\emptyset$, forcing pseudocycles and
bordisms to be absolute, we omit this word.
\interpar

If $(W,F)$ is a relative $\Ps{k}$-bordism between $(V_{0},f_{0})$ and
$(V_{1},f_{1})$, then the restriction $$(\overline{bW\setminus(V_{0}\sqcup
V_{1})}, F|\overline{bW\setminus(V_{0}\sqcup V_{1})})$$ is a
$\Ps{k-1}$-bordism between  $(bV_{0},f_{0}|bV_{0})$ and
$(bV_{1},f_{1}|bV_{1})$.
\interpar

We now form bordism
classes of relative \Ps{k}-cycles in $(X,A)$ and bordism classes of
\Ps{k}-cycles in $X$ and in $A$.  The usual argument that bordism
gives an equivalence relation works because of Corollary 6b, in \S2,
together with the techniques of gluing and smoothing maps described in the
Appendix (to verify transitivity).  Thus, we obtain bordism or pseudohomology groups
$\Ps{k}(X,A), \Ps{k}(X)$, and $\Ps{k}(A)$. As usual, by restricting a
relative pseudocycle to the boundary, we obtain a homomorphism
$\Ps{k}(X,A)\rightarrow \Ps{k-1}(A)$, which we denote by $\partial$.

\interpar

When $V$ (resp.,$W$) is compact in the definitions above we call the
pseudocycle (resp., bordism) compact. If we restrict to compact
pseudocycles and bordisms, then we obtain the classical (oriented) bordism
groups $\Omega_{k}(X,A)$, which may be defined by smooth maps, since our targets
$(X,A)$ are smooth.  Thus, we have  natural ``inclusion'' homomorphisms
$\Omega_{k}(X,A)\rightarrow \Ps{k}(X,A)$, etc.

\interpar

After the following preliminary lemma, we present some of the basic
properties of pseudohomology.
\interpar 

\begin{lemma}

Let $(X,A)$ and $(Y,B)$ be smooth manifold pairs, and let
$h:(X,A)\rightarrow (Y,B)$ be a proper, smooth map.  If
$(V,f)$ is a relative \Ps{k}-cycle in (X,A), then $(V,hf)$ is a relative
\Ps{k}-cycle in
$(Y,B)$. If
$(W,g)$ is a relative \Ps{k}-bordism in $(X,A)$, then $(W,hg)$ is a
relative \Ps{k}-bordism in $(Y,B)$. 
\end{lemma}
\noindent{\bf Proof:} Clearly $(V,hf)$ represents a smooth, relative
bordism class in $(Y,B)$ and $(W,hg)$ represents a smooth relative bordism.
The required dimension conditions follow immediately from Corollary 6e
 (\S2).
\interpar

Thus, for proper $h$, the rule $(V,f)\mapsto
(V,hf)$ defines a homomorphism \[\Ps{k}(h):\Ps{k}(X,A)\rightarrow
\Ps{k}(Y,B).\]
As one example, given the smooth manifold pair $(X,A)$ in $\mathcal{C}^{prop}$, the
inclusion
$i:A\hookrightarrow X$ gives the homomorphism $\Ps{k}(i):\Ps{k}(A)\rightarrow
\Ps{k}(X)$.  Usually, to simplify notation, we replace $\Ps{k}(h)$ by
$h_{\ast}$.
 
\interpar

Note that without the properness condition on $h$, the lemma would be
false (see Examples 1 and 2 below in this section).

\interpar 

 Lemma 4.1  immediately
gives the following corollary.

\begin{prop}[functoriality]
	\Ps{k} is a functor from $\mathcal{C}^{prop}$ to the category of abelian
goups and homomorphisms.
\end{prop}
\interpar

It will be convenient from time to time to make use of the corresponding
\emph{reduced} functor $\tilde\Ps{\ast}$, which is defined in the usual
way for pairs of the form $(X,\emptyset)$.  That is, $\tilde\Ps{\ast}(X)=
ker(\Ps{\ast}(X)\rightarrow\Ps{\ast}(pt))$.  Of course, the collapse map
$X\rightarrow pt$ must be proper, i.e., $X$ must be compact.  
\interpar

\begin{prop}[homotopy invariance]
	\Ps{k} is a homotopy functor on $\mathcal{C}^{prop}$.
\end{prop}

\noindent{\bf Proof:}  Choose any relative \Ps{k}-cycle $(V,f)$ in $(X,A)$,
and let $$H:(X,A)\times[0,1]\rightarrow (Y,B)$$ be a smooth, proper
homotopy.  Then $H\circ (f\times {id_{[0,1]})}$ is a smooth map in $\mathcal{C}$
but not in $\mathcal{S}$.  So, we smooth the corners of $V\times[0,1]$ and 
adjust the map $H(f\times {id_{[0,1]})}$ as in the Appendix, obtaining a
smooth map $G:(V\times[0,1])_{\alpha}\rightarrow Y$ which is a smooth-bordism
between $(V,H_{0}f)$ and $(V,H_{1}f)$. It remains to verify that the smooth
dimension conditions of Definition 6 are satisfied.  The construction in the
Appendix insures that the limit sets of $G$ and of $G|bV\times[0,1]$ are the
same as those of
$H(f\times id)$ and $H((f|bdV)\times id)$, respectively, so we may compute with
the latter.  But, by Corollary 5 (\S 2),
$$L(H\circ (f\times id)) = H(L(f\times id))=H(L(f)\times [0,1]),$$
from which the inequality $\ld{H\circ(f\times
id)}\leq k-1$ follows immediately. Similarly for the inequality 
$\ld{H\circ((f|bV)\times id)}\leq k-2.$
This completes the proof of the proposition.

\interpar

Of course, the results above all specialize to the absolute case.

\interpar
Now let $(X,A)$ be a pair in $\mathcal{C}^{prop}$, and let $i$ and $j$ be the
(proper) inclusions
$A\hookrightarrow X$ and $X\hookrightarrow (X,A)$, respectively.  These maps,
together with restriction to the boundary already described above, induce a
sequence 
\[\cdots\stackrel{\partial}{\rightarrow}\Ps{k}(A)\stackrel{i_{\ast}}{\rightarrow}
\Ps{k}(X)\stackrel{j_{\ast}}{\rightarrow}\Ps{k}(X,A)\stackrel{\partial}{\rightarrow}
\Ps{k-1}(A)\stackrel{i_{\ast}}{\rightarrow}\cdots,\]
which we call the $\Ps{\ast}$-sequence of the pair $(X,A)$.

\interpar

\begin{prop}[exact sequence]

The $\Ps{\ast}$-sequence of the pair $(X,A)$ is exact.

\end{prop}

The proof follows closely the standard proof of exactness for classical bordism.
However, here there are two added complications.   First, we
need to check various smooth dimension conditions to make sure the constructions
remain within the class of pseudocycles or bordisms of pseudocycles.  These
conditions are not hard to check with the aid of Proposition 3 and Corollary 6
of \S 2. Details for this are left to the reader. Secondly, in one part of the
proof of exactness at $\Ps{k}(X,A)$, we must glue together along their boundaries
a  relative \Ps{k}-cycle $(U,f)$ in $(X,A)$ and a \Ps{k-1}-nullbordism
$(V,g)$ of $(bU,f|bU)$ in $A$, smoothing both $U\cup V$ and $f\cup g$.  Calling
the result $(W,h)$, we construct $(W\times [0,1], h\circ pr)$, which is to
serve as a relative bordism between $(W,h)$ and the original $(V,f)$.  However,
first corners must be rounded and $h\circ pr$ correspondingly
smoothed, and this all must be done without violating the requisite
limit-dimension conditions. These steps all follow from results and techniques of  the
Appendix,  to which we refer the reader.  

\interpar

\begin{prop}[excision]

Let $(X,A)$ be a smooth manifold pair in $\mathcal{C}^{prop}$, and
let
$U$ be an open subset of
$X$ such that $\overline{U}\subseteq int(A)$ and such that the inclusion
$inc:(X\setminus U,A\setminus U)\hookrightarrow (X,A)$ belongs to
$\mathcal{C}^{prop}$.  Then
\[\Ps{\ast}(inc):\Ps{\ast}(X\setminus U, A\setminus U)\rightarrow\Ps{\ast}(X,A)\]
is an isomorphism.
\end{prop}
\interpar

The proof makes use of the following general fact about smooth manifold pairs. 

\begin{lemma}
Let $(X,A)$ be a smooth manifold pair such that $dimA=dimX$, and let $K$ be a
closed subset of the topological interior $A^{\circ}$ of $A$ in $X$.  Then
there exists a smooth manifold pair $(X,B)$ such that $K\subseteq
B^{\circ}\subseteq B\subseteq A^{\circ}$ and $(X,B)$ is a (strong)
smooth, proper deformation retract of $(X,A)$.
\end{lemma} 
{\bf Proof:}  This fact follows quickly from results and methods of J. Cerf ~\cite{c},
in particular, Chapitre II, \S\S2,3. We refer the reader to that paper.
\interpar

\noindent{\bf Proof of Proposition 11:}  (1)
\emph{Surjectivity}. Setting $K=\overline{U}$,  we make use of the smooth manifold pair
$(X,B)$ of Lemma 4.2, which is a smooth, proper deformation retract of
$(X,A)$ and satisfies $\overline{U}\subseteq B^{\circ}\subseteq B\subseteq
A^{\circ}$.  It follows from the homotopy property of
$\Ps{\ast}$ that, given a class
$z\in$\Ps{k}(X,A), we may represent it by a \Ps{k}-cycle
$(V,f)$ such that $f(bV)\subseteq B$.  By definition, there exist smooth maps
$F:M\rightarrow X$ and
$G:N\rightarrow A$ such that $dimM=k-2,\: dimN=k-3,\:F(M)\supseteq L(f)$, and
$G(N)\supseteq L(f|bV)$.  

We now use a theorem of Whitney (e.g., see ~\cite{br}) to construct a smooth function
$\lambda:X\rightarrow[0,1]$ such that $\lambda(B)=1$ and
$\lambda(X\setminus A^{\circ})=0$. Whitney's theorem and proof hold just as well in the presence of corners as in their absence. By Sard's Theorem, there exists a $c\in(0,1)$
that is a regular value simultaneously for the maps $\lambda f,\: \lambda
(f|bV),\:
\lambda F,\:$ and $\lambda G$.  Let $W=(\lambda f)^{-1}([0,c])$, and set $g=f|W$.
We claim that $(W,g)$ is a \Ps{k}-cycle in $(X\setminus U,A\setminus U)$ which
is \Ps{k}-bordant to $(V,f)$ in $(X,A)$. 

To verify the claim, note first that since $\lambda f(W)\subseteq[0,c]$,
$W$ is a smooth submanifold of $intV$ with boundary
$(\lambda f)^{-1}(c)$ and no corners. Further,
$g$ maps $(W,bW)$ smoothly into $(X\setminus U,A\setminus U)$, and
$(V\times[0,1], f\circ pr)$ is a smooth relative bordism in $(X,A)$ between
$(V,f)$ and
$(W,g)$. It remains only to verify the requisite dimension conditions on limit
sets.

First we verify the dimension conditions for $(W,g)$. Since $W$ is a closed
subset of
$V$, we have
$L(g)\subseteq L(f)$ (Proposition 3, \S 2), and so $\ld{g}\leq\ld{f}\leq k-2,$
as required. Further, since
$bW=\lambda f^{-1}(c)$,  we have $\ld{g|bW}=\ld{g|(\lambda
f)^{-1}(c)}=\ld{f|f^{-1}(\lambda^{-1}(c))}$.  By construction, both $f$ and
$F$ are transverse to $\lambda^{-1}(c)$. Thus, applying Corollary 7 of \S 2, we
get $\ld{g|bdW}\leq dimM-1\leq k-3,$ as required. 

The dimension conditions for $(V\times[0,1], f\circ pr)$ are verified more
easily. Thus, $\ld{f\circ pr}=\ld{f}\leq k-2 < k-1,$ (Corollary 6, \S 2), as
required. Also, again using Corollary 6,
\begin{eqnarray*}
\ld{f\circ
pr|\overline{b(V\times[0,1])\setminus((V\times 0)\sqcup (W\times
1))}} & = & \ld{f\circ pr|(bV\times[0,1])\cup\overline{V\setminus
W}\times 1} \\
 & = & \max\{\ld{f|b V}, \ld{f|\overline{V\setminus W}}\}\\
 & \leq & k-2.
\end{eqnarray*}
This completes the proof of surjectivity.

\interpar

(2) \emph{Injectivity}. Let $B$ and $\lambda$ be as in the proof of
surjectivity.  The smooth, proper deformation retraction
$(X,A)\rightarrow (X,B)$ of that proof restricts to a smooth proper
deformation retraction
$(X\setminus U,A\setminus U)\rightarrow (X\setminus U, B\setminus U)$.  It
follows that any class
$w\in\Ps{k}(X\setminus U,A\setminus U)$ that maps to $0\in\Ps{k}(X,A)$ can be
represented by a \Ps{k}-cycle in $(X\setminus U,B\setminus U)$ that is
\Ps{k}-nullbordant in
$(X,B)$. Accordingly, let
$(V,f)$ be such a
\Ps{k}-cycle in
$(X\setminus U, B\setminus U)$ representing $w$, and let $(W,g)$ a
\Ps{k}-nullbordism of
$(V,f)$ in
$(X,B)$.  Recall that we have $\lambda(B)=1$ and
$\lambda(X\setminus A^{\circ})=0$.

Now we again make use of smooth maps provided by the definitions,
$F:M\rightarrow X\setminus U,\: F\pri:M\pri\rightarrow A\setminus U,\:
G:N\rightarrow X,\mbox{and}\: G\pri:N\pri\rightarrow A,$ which satisfy the
following:
$dimM=k-2\:
\mbox{and}\:F(M)\supseteq L(f);\: dimM\pri=k-3\:\mbox{and}\:
F\pri(M\pri)\supseteq L(f|bV);\: dimN=k-1\:\mbox{and}\: G(N)\supseteq
L(g);\: dimN\pri =k-2\:\mbox{and}\: G\pri(N\pri)\supseteq
L(g|\overline{bW\setminus V}).$  Using the function $\lambda$, we
choose a $d\in (0,1)$ that is simultaneously a regular value for
$\lambda f,\: \lambda f|bV,\: \lambda g,\:  \lambda g|bW,\: 
\lambda F,\: 
\lambda F\pri,\:  \lambda G,\: \mbox{and}\:  \lambda G\pri.$ 

Note that
$(\lambda g)^{-1}(d)$ separates $W$(resp.,$bW$) into two submanifolds
$W_{\pm}$ (resp., $(bW)_{\pm}$), with the notation chosen so that
$W_{+}=(\lambda g)^{-1}([0,d]),$ etc.  Thus, $(bW)_{+}\subseteq
int(V)$.  Set
$h=g|W_{+}$.
 The remainder of the proof now consists of verifying the following two claims:
\begin{enumerate}
 \item \emph{$((bW)_{+}, h|(bW)_{+})$ is \Ps{k}-bordant to
$(V,f)$ in
$(X\setminus U,A\setminus U)$}.
 \item \emph{$(W_{+},h)$ is a \Ps{k}-nullbordism of $((bW)_{+},
h|(bW)_{+})$ in $(X\setminus U,A\setminus U)$}.
\end{enumerate}

For the first claim, we make use of the smooth bordism $(V\times[0,1],f\circ
pr)$ between $(V,f)$ and $((bW)_{+}, h|(\partial W)_{+})$ in
$(X\setminus U,A\setminus U)$. Of course, corners must be straightened and the
map smoothed so as to preserve limit dimension.  Again, this is done by the
techniques of Appendix A.  The verification of dimension conditions proceeds
similarly to that in the proof of surjectivity.  The only slightly subtle point
involves
$\ld{h|f^{-1}(d)}$, which is evaluated with the aid of Corollary 7, \S 2, as
before.  Similar arguments apply to the second claim.

This concludes the proof of the excision property.
\interpar

\begin{corr}
$\Ps{\ast}$ takes finite disjoint unions to finite direct sums.
\end{corr}

The proof is a well-known induction on the number of terms in the union, which works for
any functor satisfying the Eilenberg-Steenrod axioms.  For the induction step, one looks
at the long exact sequence for the pair (full union, union with one term deleted). The
excision property allows one to remove the subspace, etc.
\interpar

Suppose now that $\{X_{i}\}$ is a handle decomposition for the pair $(X,A)$ (cf., \S 3).  
Combining the foregoing results, we compute just as for homology:
$$\Ps{k}(X_{i},X_{i-1})\approx
\bigoplus_{j}\Ps{k}(D^{i}_{j},S^{i-1}_{j}),$$ where the subscript $j$ indexes the
$i$-handles attached to $X_{i-1}$.  To complete this calculation, we need to
compute the groups $\Ps{k}(D^{i},S^{i-1})$. 
This can be done with the aid of some of the computational facts that we now present. Let
$\iota:\Omega_{\ast}(X,A)\rightarrow\Ps{\ast}(X,A)$ denote the homomorphism
induced by the inclusion of bordism classes.  There is also a corresponding
homomorphism $\tilde\iota$ of the reduced theories.  Recall that there is a natural
transformation of functors $\mu:\Omega_{\ast}\rightarrow H_{\ast}$ defined by evaluating
fundamental classes.
\interpar

\begin{prop}\label{prop:point}
\begin{enumerate}
  \item If X is compact, then $\iota:\Omega_{0}(X)\rightarrow
\Ps{0}(X)$ is an isomorphism.  Therefore, the composite
$\iota\circ(\mu)^{-1}:\Hm{0}(X)\rightarrow\Ps{0}(X)$ is an isomorphism.
  \item If $X$ is $n$-dimensional, then
$\Ps{k}(X)=0$, for all
$k>n$.
\end{enumerate}

\end{prop}
{\bf Proof:}  1. If $(V,f)$ (resp., $(W,g)$) is a \Ps{0}-cycle (resp.,
\Ps{0}-bordism) in $X$, then $V$ (resp., $W$) is compact (Proposition 3a). 
It follows that
$\iota$ is just the identity homomorphism.

2. By definition, when $k>n$, every smooth bordism in $X$ is a \Ps{k}-bordism.
Therefore, given any \Ps{k}-cycle $(V,f)$ in $X$, we can find a \Ps{k}-nullbordism of
$(V,f)$: namely,
$(V\times [0,\infty), f\circ pr)$.

This completes the proof.

\interpar

In the special case $X=pt$, Proposition 13 implies the dimension axiom asserted
in Theorem 1.  Thus, combining this proposition with Propositions 8---11, we
have completed the proof of Theorem 1.

\interpar 

Both assertions of Proposition 13 carry over to the reduced groups. Note
that the first assertion immediately implies the following:  \emph{if
$f:X\rightarrow Y$ establishes a bijection between the components of the compact
manifolds $X$ and $Y$, then $f_{\ast}:\Ps{0}(X)\rightarrow\Ps{0}(Y)$ is an isomorphism.}
\interpar

\begin{prop}\label{prop:ball-point}
 \begin{enumerate} 

 \item  The inclusion $pt=D^{0}\hookrightarrow D^{n}$ induces an isomorphism
$\Ps{\ast}(pt)\approx\Ps{\ast}(D^{n})$.
 \item  There is an isomorphism $\Ps{\ast}(D^{n},S^{n-1})\approx
\tilde{\Psi}_{\ast-n}(S^{0})$, which is a composition of boundary maps of long
exact sequences and excision isomorphisms.
   
\end{enumerate}
\end{prop}

The proof is the same as for ordinary homology.
\interpar

The remaining computational observations imply items a and b, listed in the
Remark in \S1.
\interpar

\noindent{\bf Example 1:}  Let $U=V$ denote the open unit interval $(0,1)$, and 
define $f:U\rightarrow V$
and $g:V\rightarrow S^{1}$ by $f(x)=x$ and $g(y)=\exp(2\pi iy)$. Then
$L(f)=\emptyset$, so $(U,f)$ defines a $\Ps{0}$-cycle in $V$.  On the other hand,
$L(gf)\not=\emptyset$, so $(U,gf)$ is not a $\Ps{0}$-cycle in $S^{1}$. 
\interpar

This example shows that composition with maps that are not proper does not
define a homomorphism of pseudohomology.
\interpar

\noindent{\bf Example 2:}  Let $U=W=S^{1},$ the unit sphere, and let
$V=U\times R$.  Let $i$ denote the inclusion map $U\hookrightarrow V$ given by
$x\mapsto(x,0)$, and let $p:V\rightarrow W$ be the projection. As we shall see
below, the pair $(U,id_{U})$ represents a non-zero class in \Ps{1}(U), and so
$(p\circ i)_{\ast}[U,id_{U}]\not=0$ in $\Ps{1}(W)$.  However, $(U,i)$ is
\Ps{0}-nullbordant in $V$ (with null bordism the inclusion
$S^{1}\times[0,\infty)\hookrightarrow S^{1}\times R$), so that,if defined,
the homomorphism $p_{\ast}$ satisfies $p_{\ast}i_{\ast}[U,id_{U}]=0$.
\interpar

Therefore, in the absence of properness, \Ps{\ast} does not define a functor.

\begin{prop}\protect\label{prop-noninject}

Let $X$ be any smooth $k$-manifold with empty boundary, let $V=X$, and let
$f:V\rightarrow X$ be any diffeomorphism.  Then $(V,f)$ is a
\Ps{k}-cycle in $X$ that does not represent $0\in \Ps{k}(X)$. 
\end{prop}
Thus, $\Ps{k}(X)\not=0$, whereas $H_{k}(X)=0$ .
It follows that when $X$ has non-compact components,
there is no natural surjection
$\Hm{k}(X)\rightarrow \Ps{k}(X)$.
\interpar

\noindent{\bf Proof of Proposition 15:}  Since $f$ is a
diffeomorphism,
$L(f)=\emptyset$, so $(V,f)$ is a \Ps{k}-cycle. The remainder of the proof is
similar to Hirsch's proof of the Brower Fixed Point Theorem, as presented in ~\cite{m2}.  Suppose then that, counter to our proposed conclusion, we have 
 a
\Ps{k}-nullbordism of
$(V,f)$ in $X$, say $(W,g)$. Since $g|V=f$,  $g$ is surjective.  We observe that
$sdL(g)\leq k-1$, so that
$L(g)$ has Lebesgue measure $0$ in $X$.  Similarly, by Sard's Theorem, the set of
critical values of $g$ has measure $0$.  It follows that there exists a regular
value $p$ of $g$,  that is not in $L(g)$. By the definition of $g$, $p$ is also a
regular value of
$g|bW=g|V=f$. Set $g^{-1}(p)=A$. Then $A$ is a non-empty, compact
(Proposition 3a, \S 2),
$1$-dimensional submanifold of
$W$ such that $bA = A\cap bW=\{f^{-1}(p)\}$, a singleton subset of $A$. However, there is
no such manifold, because every compact $1$-manifold has an
\emph{even} number of boundary points.  This contradiction proves the claim.
\interpar

\begin{prop}
Let $X$ be any smooth manifold, and let $j:X\rightarrow
X\times R$ denote the (proper) inclusion given by $x\mapsto(x,0)$. Then
$j_{\ast}:\Ps{\ast}(X)\rightarrow\Ps{\ast}(X\times R)$ is trivial.
\end{prop}
{\bf Proof:}  If $(V,f)$ is a \Ps{k}-cycle in $X$, then $(V\times [0,\infty), f\times
inc)$ is a \Ps{k}-nullbordism of $(V,f)$ in $X\times R$.

\begin{corr}
Suppose that $X$ is a compact, orientable $n$-manifold with
empty boundary, and suppose that 
$\psi:H_{\ast}\rightarrow\Ps{\ast}$ is a natural tranformation of functors on
$\mathcal{C}^{prop}$. Then $\psi:H_{\ast}(X\times R)\rightarrow
\Ps{\ast}(X\times R)$ is non-injective.
\end{corr}
{\bf Proof:} The homomorphism $j_{\ast}:H_{n}(X)\rightarrow H_{n}(X\times R)$
is a non-zero isomorphism.  The result is now immediate from the foregoing
proposition.

\section{Smoothing circuits}
We describe the representation of homology by oriented, singular
$k$-circuits due to Rourke and Sanderson~\cite{rs} (see also ~\cite{brs}) and then pass
to weakly piecewise smooth (WPS) circuits .  These can be used to
represent the homology of smooth manifold pairs, and they have the virtue
of allowing us some control over smooth dimension.  We conclude by showing
how to smooth WPS circuits so as to produce pseudocycles.

\subsection{Circuits}

We begin by reminding the reader of some basic terminology from
piecewise linear (PL) topology.

\begin{defn}
All of the simplicial complexes that we consider will be locally finite
collections of simplexes contained in some finite dimensional Euclidean
space.  If $K$ is a simplicial complex and $S$ is a subset of $K$, then
we let $|S|$ denote the union of all open simplexes $\sint{\sigma}$ for
$\sigma\in S$.  We call $|K|$ a \emph{polyhedron} and $K$ a
\emph{triangulation} of
$|K|$.  The set $St(S,K)$ consists of all $\sigma\in K$ having some face
in $S$.  The subset $|St(S,K)|\subseteq|K|$ is always open; we sometimes
call it the (open) star of $|S|$ in $|K|$.

We let $S(i)$ denote the collection of $i$-simplexes in $S$, and $S^{i}$
the union of all $S(j),\,j\leq i$. If $S$ is a subcomplex of $K$, then
so is $S^{i}$. 

A polyhedral pair $(Q,P)$ is a pair of topological spaces for which there
exist a simplicial complex $L$ and subcomplex $K$ with $Q=|L|$ and
$P=|K|$.  The pair $(L,K)$ is said to triangulate $(Q,P)$.  The dimension
of $(Q,P)$ is the usual dimension of $L$.  This is independent of the
choice of triangulating $L$.
\end{defn}

\begin{defn}
\begin{enumerate}
  \item A \emph{relative $k$-circuit} is a \emph{compact} polyhedral pair
$(Q,\delta Q)$ such that there is a $(k-2)$-dimensional subpolyhedron
$S(Q)\subseteq Q$ satisfying:
    \begin{enumerate} 

      \item $Q= \overline{Q\setminus S(Q)}$.

      \item $Q\setminus S(Q)$ is a piecewise-linear $k$-manifold
with boundary $\delta Q\setminus S(Q)$.

      \item $\delta Q$ is a relative $(k-1)$-circuit with
$\delta(\delta Q)=\emptyset$ and $S(\delta Q)=S(Q)\cap \delta Q.$

If $\delta Q=\emptyset$, we say that $Q$ is a closed or
absolute $k$-circuit (or simply a $k$-circuit).  Then (c) may be
rephrased as

      \item  $\delta Q$ is a $(k-1)$-circuit with $S(\delta
Q)=S(Q)\cap\delta Q$.
    \end{enumerate}

  \item A \emph{null
bordism} of a relative $k$-circuit $Q$ is a relative $(k+1)$-circuit
$R$ such that $Q\subset\delta R$ and $S(R)\cap Q=S(Q).$
It follows that $S(R)\cap \delta Q=S(\delta R)\cap \delta
Q=S(\delta Q)$ and that $Q\setminus S(Q)$ and $\delta Q\setminus S(\delta
Q)$ are properly embedded PL submanifolds of $\delta R\setminus S(\delta
R)$. Absolute bordisms are the obvious specializations.  For example, if
$R$ is a null bordism of the relative $k$-circuit $Q$, then $\delta
R\setminus (Q\setminus \delta Q)$ is a null bordism of the absolute
circuit $\delta Q$.

  \item A \emph{singular $k$-circuit in a space $Z$} is a pair
$(P,a)$,  where P is a $k$-circuit and $a:P\rightarrow Z$ is a continuous
map.  Similarly for bordisms of singular $k$-circuits. The same
terminology applies to the relative case.

  \item We say that a circuit or bordism of circuits $T$ is
oriented if the manifold
$T\setminus S(T)$ is oriented; then it induces an orientation on the
manifold $\delta T\setminus S(T)$, so that $\delta T$ inherits an
orientation from $T$, etc. 
\emph{In this paper we restrict attention entirely to oriented circuits
and bordisms although we may not mention this further explicitly.}

\end{enumerate}
\end{defn}

With these definitions, one can now define groups of bordism classes of
relative (and absolute) $k$-circuits.  If $(X,A)$ is a pair of spaces,
then we denote its circuit bordism group by $\mathcal{H}_{k}(X,A)$, and
similarly in the absolute case.

\interpar

Let $(Q,a)$ be a singular, relative $k$-circuit in $(X,A)$.  As
stated in Definition 9d, $Q\setminus\delta Q$ is oriented.  It
follows that $Q$ carries a fundamental (singular) homology class. (To see
this in the absolute case, triangulate $Q$ such that $S(Q)$ is a
subcomplex, and orient every
$k$-simplex compatibly with $Q$.  The formal sum of these is a simplicial
$k$-chain on $Q$, which one quickly computes to be a $k$-cycle.  This
cycle can now be used to define a singular $k$-cycle that represents
$[Q]$. The relative case is similar.)  The class $[Q]$ allows us to
define an evaluation map $e:\mathcal{H}_{k}(X,A)\rightarrow H_{k}(X,A)$
by the usual rule $e[Q,a]=a_{\ast}([Q])$.
\interpar

\noindent{\bf Theorem (Rourke and Sanderson~\cite{rs}):} \emph{The evaluation map 
$e:\mathcal{H}_{k}(X,A)\rightarrow H_{k}(X,A)$ is a natural isomorphism for all k.}
\interpar

This theorem provides us with almost the right context in which to
compare homology with pseudohomology.  However, notice that circuit
homology uses continuous maps, whereas pseudohomology uses smooth
maps with constraints on their limit dimensions. Therefore, it is
convenient to modify circuit homology slightly to adapt to
this.  Our first step is to introduce the notion of a
\emph{weakly piecewise-smooth} map (which we henceforth call a WPS
map).

\begin{defn}

Let $S$ be a subset
of a locally-finite simplicial complex $K$,
let $Y$ be smooth, and let  $f:|S|\rightarrow Y$ be a continuous map.  We
say that $f$ is weakly piecewise-smooth with respect to $S$ (WPS
w.r.t.$S$) if $f|\sint{\sigma}:\sint{\sigma}\rightarrow Y$ is smooth
for every simplex $\sigma$ in $S$.\footnote{The usual notion of
piecewise smooth map, which we \emph{do not} want here, requires that $f$
be smooth on each $\sigma$.}  

If $P$ is a polyhedron,
then we say that $f$ is WPS if $f$ is WPS w.r.t. $K$ for some $K$
triangulating $P$. 

\end{defn}

The following properties of this definition are easily checked: 
\begin{itemize}
  \item If $f$ is
WPS w.r.t.
$K$,  and
$K^{\prime}$ is a subdivision of $K$, then $f$ is WPS
w.r.t. $K^{\prime}$.  
  \item If $f:P\rightarrow Y$ is WPS, then for every $L$
triangulating
$P$, there is a subdivision $L^{\prime}$ of $L$ such that $f$ is
WPS w.r.t. $L^{\prime}$.
\end{itemize}

\interpar

At this point, we observe that we can repeat Definition 9c, with the
topological space $Z$ replaced by a smooth manifold and continuous maps
replaced by WPS maps; similarly for the relative case.  This results in
bordism groups
$\mathcal{H}_{k}^{WPS}(X,A)$ and an ``inclusion'' homomorphism
$$\jmath:\mathcal{H}_{k}^{WPS}(X,A)\rightarrow
\mathcal{H}_{k}(X,A),$$
for each smooth manifold pair $(X,A)$.

\begin{prop}

The homomorphism $\jmath$ is an isomorphism. 

\end{prop}

\noindent {\bf Proof (sketch):}  For simplicity, start with the
absolute case. An induction argument reduces the proof to showing that if
$(Q,P)$ is a compact, polyhedral pair and $f:Q\rightarrow X$ a continuous
map such that $f|P$ is WPS, then there is a homotopy $relP$ (i.e.,
pointwise stationary on $P$) between $f$ and a WPS map. In turn,
another induction argument reduces this to the special case in which $Q$
is a simplex $\sigma$ and $P$ its boundary
$\dot{\sigma}$. So suppose that we are in this situation, choose a
topological metric on
$X$, and then choose a smooth $\epsilon$-approximation to
$f|\sint{\sigma}$, say
$g^{\prime}$. Here $\epsilon$ is a continuous function on $\sigma$
which is positive on $\sint{\sigma}$ and $0$ on 
$\dot{\sigma}$. We extend
$g^{\prime}$  to a  map $g$ on $\sigma$ by setting $g=f$
on
$\dot{\sigma}$, and we check that $g$ is WPS.
Moreover, by choosing
$\epsilon$ suitably small, it is easy to arrange that $g$ be homotopic to
$f\: rel\dot{\sigma}$. This concludes the sketch in the absolute case. 

The relative case can be derived from the absolute case.  We start with
$Q, P,$ and
$f$ as before, only now we have a smooth manifold pair $(X,A)$ and
subpolyhedra
$Q^{\prime}\subseteq Q$ and $P^{\prime} = P\cap Q^{\prime}$, with $f$ a
map of pairs
$(Q,Q^{\prime})\rightarrow (X,A)$. We want to conclude that there is a
homotopy $rel\: P$ of maps of pairs between $f$ and a WPS map. To obtain
this, first apply the absolute case to $Q^{\prime}, P^{\prime}$, and
$f|Q^{\prime}$.  Then apply the homotopy extension property to obtain
a homotopy of pairs $rel\: P$ between $f$ and a continuous map which is
WPS on
$P\cup Q^{\prime}$.  Then apply the absolute case again to the pair $(Q,
P\cup Q^{\prime})$.

This concludes our sketch.

\interpar

It follows that we may use WPS circuits to represent the homology of
smooth manifold pairs.

\interpar

We conclude this subsection with a result which shows that smooth
dimension can be defined in terms of WPS maps.

\begin{prop}

Let $Y$ be smooth, and let $S$ be a subset of $Y$.  Then $sdS$ is the
smallest integer $m$ for which $S$ is contained in the image of a WPS map
$P\rightarrow Y$ with $dimP\leq m$.

\end{prop}

\noindent{\bf Proof:}  The proposition is clearly equivalent to
the statement,
\begin{quote}
		$sdS\leq m$ iff $S$ is in the image of a WPS map $P\rightarrow Y$
with $dimP\leq m$,
\end{quote}
for which we now give a proof. If $sdS\leq m$, then $S$ is in
the image of a smooth map $X\rightarrow Y$ with $dimX\leq m.$  Use a
smooth triangulation of $X$ to transform this into a WPS map.  On the
other hand, if there is a WPS map $f:P\rightarrow Y$ with image
containing $S$ and with $dimP\leq m$, form a smooth map as follows. Let
$dimP=k\leq m$, and choose a triangulation $K$ of $P$ such that $f$ is WPS w.r.t. $K$.
Choose any
$\sigma\in K(i)$. The smooth $k$-manifold
$\sint{\sigma} \times R^{k-i}$ maps smoothly to $Y$ by first projecting to
$\sint{\sigma}$ in
$P$ and then mapping to $Y$ by $f$. The disjoint union of
all these maps is a smooth map to $Y$ from a $k$-manifold 
with image containing $S$.  Thus, $sdS\leq m$.  This completes
the proof.

\subsection{Manifolds in circuits}

The Rourke-Sanderson definition of a $k$-circuit $Q$ is independent of
the choice of simplicial complex $L$ triangulating $Q$. Thus, for
example, the singular subset $S(Q)\subseteq Q$ is simply defined to be a
$(k-2)$-dimension subpolyhedron of $Q$ and \emph{not} necessarily the
$(k-2)$-skeleton of a triangulation.  While this invariance is undoubtedly
a useful feature of their theory, it poses some obstacles for us from the
viewpoint of smoothing the complement of the singular set. We deal with
these by choosing arbitrary polyhedral triangulations of $Q$ (i.e., not necessarily
respecting $S(Q)$) and choosing $(k-2)$-dimensional subcomplexes that both
contain the ``bad points'' in $S(Q)$ and have smoothable complements. 
By ``bad points'' we mean non-manifold points, which we now describe.

\begin{defn}
Let $K$ be a $k$-dimensional simplicial
complex. We call a point $p$ in $|K|$ a \emph{manifold-point} if it has a
neighborhood in $|K|$\,that is PL-homeomorphic to $R^{k}$ or $R^{k}_{+}$
(closed, Euclidean half-space). Otherwise we say that $p$ is a
\emph{non-manifold point}.
\end{defn}

Clearly the non-manifold points form a
closed subset of $|K|$, but we can say more than this.  The
following lemma presents a standard fact from PL topology; we
provide a proof for the reader's convenience.

\begin{lemma}

The non-manifold points of $|K|$ form a subpolyhedron triangulated by
a subcomplex of $K$.

\end{lemma}

\noindent{\bf Proof:}  If $p,q\in |K|$, both belonging to the same
open simplex, then it is easy to construct a PL homeomorphism
$f:|K|\rightarrow |K|$ such that $f(p)=q$.  It
follows that both $p$ and $q$ are manifold points or non-manifold
points.  Therefore, the non-manifold points of
$|K|$ consist of a union of open simplexes.  Since this set is closed,
it must be a union of closed simplexes, completing the proof.

\interpar

In the remainder of this subsection, we shall be consistently dealing with
the following separate three cases:
\begin{enumerate}
\item $P$ is an absolute $(k-1)$-circuit triangulated by a finite
simplicial complex $K$.
\item $(Q, P)$ is a relative $k$-circuit triangulated by $(L,K)$.
\item $R$ is a relative $(k+1)$-nullbordism of $(Q,P)$ with $(R,\delta R)$
triangulated by $(N,M)$ and $L$ and $K$ subcomplexes that triangulate
$Q$ and $P$ respectively. 
\end{enumerate}

In each of these cases, we define a subpolyhedron $\Sigma$ of codimension
two, and we show that it contains all the bad points.
\interpar

\noindent{\bf Definition of $\Sigma$:}  
\begin{enumerate}
  \item $\Sigma=|K^{k-3}|$.
  \item $\Sigma=|L^{k-2}\setminus K(k-2)|$.
	 \item $\Sigma=|N^{k-1}\setminus(M(k-1)\cup St(K(k-2),N))|$.
\end{enumerate}
The reader can easily convince herself/himself that $\Sigma$ is a compact
subpolyhedron of codimension two.
\interpar

\begin{prop}
\begin{enumerate}
  \item $P\setminus\Sigma$ is a PL $(k-1)$-manifold with empty boundary.
  \item $Q\setminus\Sigma$ is a PL $k$-manifold with boundary equal to
$P\setminus\Sigma$.
	 \item $R\setminus\Sigma$ is a PL $(k+1)$-manifold with boundary
$\delta R\setminus\Sigma$.  Moreover, $Q\setminus\Sigma$ is a properly
embedded, codimension-$0$ submanifold of $\delta R\setminus\Sigma$
with boundary $P\setminus\Sigma$.
\end{enumerate}
\end{prop}
{\bf Proof:}  We show first that each of the stated complements of
$\Sigma$ is a PL manifold.  Since each complement is an open set, it
suffices to show that it contains only manifold points.

By definition, in each case the non-manifold points are contained in
a codimension-two subpolyhedron ($S(P),\, S(Q),\, S(R)$, as the case
may be).  Since they form a subcomplex, they must be contained in
$K^{k-3}$ (Case a), or $L^{k-2}$ (Case b), or $N^{k-1}$ (Case c). 
Already this gives the desired result in Case a.  In Case b, the only
potentially bad points that $\Sigma$ misses are in $|K^{k-2}|\subseteq
P$.  However, since $S(Q)\cap P=S(P)$, which has codimension two in
$P$, the points in $|K^{k-2}|$ are all manifold points, disposing of
Case b.  This argument also shows that in Case c, all points except
possibly those in $St(K(k-2),N)$ are manifold points, so it remains to
deal with these.  First note that all points in $K(k-2)$ are
manifold points, since $S(R)$ meets $P$ in a codimension two
subpolyhedron.  Next, consider any $\sigma\in St(K(k-2),N)$.  We claim
that $\sint{\sigma}$ consists entirely of manifold points of $R$.  For
by Lemma 5.1, $\sint{\sigma}$ is either disjoint from or contained in
the set of non-manifold points.  In the latter case, each of its
$(k-2)$-faces would also consist of non-manifold points, which we just
showed does not happen.  This verifies the claim and concludes the
proof that all the complements of $\Sigma$ are PL manifolds.  It
remains to verify the remaining statements involving boundaries of
these manifolds.

\emph{Case a}.  Choose any $p\in P\setminus\Sigma$, and let $\sigma\in K$
be its carrier, i.e., $p\in\sint{\sigma}$.  Since $dim\sigma\geq k-2$,
there is a $q\in\sint{\sigma}$ that is not in $S(P)$, i.e., $q$ is a
non-boundary point in the PL manifold $P\setminus S(P)$.  Since, as in
the argument for Lemma 5.1, there is a PL homeomorphism $P\rightarrow
P$ that takes $p$ to $q$, $p$ cannot be a boundary point of
$P\setminus\Sigma$. 

\emph{Case b}. If
$p$ belongs to $Q\setminus\Sigma$ but not to $P$, then it belongs to an
open simplex
$\sint\sigma$ of dimension $k$ or $(k-1)$ that does not meet
$P$.  As before, we find a PL homeomorphism $f$ of $Q$ such that
$f(p)$ is in $\sint{\sigma}$ and not in $S(Q)$.  Therefore, $f(p)$ is in
the manifold interior of $Q\setminus S(Q)$, and so $p$ is in the manifold
interior of $Q\setminus\Sigma$.  This shows that the
boundary of
$Q\setminus\Sigma$ is contained in $P\setminus\Sigma$. For the reverse
inclusion, choose any $p\in P\setminus\Sigma$ and argue as before to
find a PL homeomorphism $f:Q\rightarrow Q$ keeping $p$ in its carrier in
$P$ and mapping it to $q\in P\setminus S(Q)$.  By hypothesis, $q$ is a
boundary point in $Q\setminus S(Q)$.  Hence, $p$ is a boundary point in
$Q\setminus\Sigma$.

\emph{Case c}.  The arguments here are similar to the foregoing and
will be left to the reader.

This completes our proof.
\interpar

\noindent{\bf Remark:}  We note for future reference that each of the
manifolds $P\setminus\Sigma$,\, $Q\setminus\Sigma$,\, $R\setminus\Sigma$
is contained in the complement of a skeleton:
\begin{itemize}
  \item $P\setminus\Sigma\subseteq |K\setminus K^{k-3}|$.
  \item $Q\setminus\Sigma\subseteq |L\setminus L^{k-3}|$.
  \item $R\setminus\Sigma\subseteq |N\setminus N^{k-3}|$.
\end{itemize}

\subsection{Smoothing circuits}

Smoothing a singular circuit has two parts: smooth the domain (away from
a codimension-two subcomplex), and then smooth the map.  Similarly with
singular circuit bordisms. In this subsection, we show how to smooth the
domains.

We continue to deal with the three cases
listed in the previous subsection and use the notation introduced there.
By a \emph{smoothing} of a PL manifold $T$ we mean a differentiable
structure on $T$, say $\alpha$, such that the identity map $T\rightarrow
T_{\alpha}$ restricts to a smooth embedding of each closed
top-dimensional simplex of a triangulation of $T$.  (Here $T_{\alpha}$
denotes the differentiable manifold with underlying space $T$ and
differentiable structure $\alpha$.) By a \emph{concordance} of smoothings
$\alpha,\beta$ \, of $T$, we mean a smoothing of $T\times R$  that
coincides with the product smoothing $\alpha\times R$ on $T\times
(-\infty,\epsilon)$ and $\beta\times R$ on $T\times
(1-\epsilon,\infty)$, for some suitable, small, positive $\epsilon$.  We
usually abbreviate this by saying that it is a smoothing of
$T\times[0,1]$ that coincides with
$\alpha$ on
$T\times 0$ and $\beta$ on $T\times 1$. Concordances of smoothings define
an equivalence relation on the smoothings of $T$.  This notion has obvious
refinements to the relative case.

\begin{prop}\protect\label{prop:smoothdomain}
Consider the $(k-1)$-circuit $P$, the relative $k$-circuit $Q$, and the
relative nullbordism $R$ of $Q$ discussed in 5.2 above, and also let
$\Sigma$ denote the singular subpolyhedra defined there. Then:
\begin{enumerate}
\item $P\setminus\Sigma$ admits a smoothing, which is unique up to
concordance.
\item Every smoothing of $P\setminus\Sigma$ extends to a smoothing of
$Q\setminus\Sigma$, which is unique up to relative concordance.
\item Every smoothing of $Q\setminus\Sigma$ extends to a smoothing of
$R\setminus\Sigma$, which is unique up to relative concordance. 
\end{enumerate}
\end{prop}
{\bf Proof:}  The proof in all three cases follows the same
pattern.\footnote{An easy direct construction gives a smoothing in Case
a (cf., ~\cite{ms}, pp.90--91).  However, existence in the other cases, as well as the
uniqueness assertions, is less transparent.  So we give auniform proof
for all cases.} Let
$(T,S)$ be one of the three pairs
$(P\setminus\Sigma,\emptyset),\, (Q\setminus\Sigma,P\setminus\Sigma),\,
(R\setminus\Sigma,Q\setminus\Sigma)$. We use the obstruction theory of
Hirsch~\cite{h} to prove the existence of the stated smoothing.  The
obstructions in this theory reside in local cohomology groups
$H^{i}(T,S;\Gamma_{i-1})$, with coefficents defined as follows:
$\Gamma_{n}$ is obtained from the group of all diffeomorphisms
$Diff(S^{n-1})$ by factoring out the subgroup of diffeomorphisms that
extend to $D^{n}$.  These groups are known to be finite and abelian for
all $n$,~\cite{km}, and moreover, \[\Gamma_{n}=0, \mbox{for all}\, n\leq 6.\]
The proof consists in showing that each pair $(T,S)$ has the homotopy
type of a CW pair of  small dimension (in particular, at most one, two,
and three in the three respective cases).  Thus the obstruction groups
all vanish identically, and so the desired smoothings exist.  

The uniqueness assertions are proved similarly.  Here obstructions lie in
the groups $H^{i}(T,S;\Gamma_{i})$, which vanish identically.

It remains, then, to check that $(T,S)$  has the stated homotopy type.
We do this via three lemmas.

\begin{lemma}

Suppose that $K$ is a finite simplicial complex of dimension $k$ and
$S$ is a subset of $K$ containing the $r$-skeleton $K^{r}$, for some $r$. 
Then
$|K\setminus S|$ has the homotopy type of a finite CW complex of
dimension $\leq k-r-1$.

\end{lemma}

\noindent{\bf Proof:}  We begin by proving the result when $S=K^{r}$.
Let  $K^{\prime}$ denote the barycentric subdivision of $K$, and, for
every $\sigma\in K$, let $\hat{\sigma}$ denote the barycenter of
$\sigma$.  For simplexes $s$ and $t$, write $s<t$ to indicate that $s$ is a proper face
of $t$.  Now, for any $\tau\in K$, let
$$\dot{D}(\tau)=\{\hat{\tau_{1}}\hat{\tau_{2}}\cdots\hat{\tau_{s}}
|\tau<\tau_{1}<\tau_{2}<\cdots<\tau_{s}\}.$$
Clearly, $\dot{D}(\tau)$ is a subcomplex of $K^{\prime}$.  If $\tau$
has dimension $d$, then every simplex of $\dot{D}(\tau)$ has
dimension $\leq k-d-1$.  Thus, setting $$L=\bigcup_{dim\tau=r}
\dot{D}(\tau),$$ $L$ is a subcomplex of $K^{\prime}$ of dimension
$\leq k-r-1$.  We show that $|K\setminus K^{r}|\simeq|L|$ (where
$``\simeq''$ denotes homotopy equivalence).

Every simplex of $K^{\prime}$ decomposes uniquely as a join
$\lambda\ast\mu$, with $\lambda\in(K^{r})^{\prime}$ and $\mu\in L$. 
(We allow  either $\lambda$ or $\mu$ to be empty, so that, for
example, if $\lambda=\emptyset$, we get $\mu\in L$, etc.)  To see
this, write the simplex uniquely as
$\hat{\tau_{1}}\cdots\hat{\tau_{s}}$, where
$\tau_{1}<\cdots<\tau_{s}$. Then $\lambda=\tau_{1}\cdots\tau_{\ell}$
and $\mu=\tau_{\ell+1}\cdots\tau_{s}$, where $\ell$ is the maximal
subscript $i$ such that $dim\tau_{i}\leq r$.  It follows that we may
deform $|K\setminus K^{r}|$ to $|L|$ by sliding points uniformly along
join lines.

We now deal with the simplexes in $S\setminus K^{r}$. Since each simplex
$\sigma$ of $K$ of dimension $>r$ is a union of join lines beginning in
$\sigma\cap|K^{r}|$ and ending in
$\sigma\cap|L|$ and meeting at most at endpoints, and each
open join line in the decompositions
$\lambda\ast\mu$ is contained in a unique open
simplex of $K$ (indeed of $K^{\prime}$), it follows that the
deformation described in the last paragraph leaves $|K\setminus S|$
invariant and deforms it to $|L\setminus S|$.  Observe that this
last is a union of open simplexes of $L$. The proof is now completed with
the following lemma.

\begin{lemma}
Let $B$ be a finite simplicial complex and $A\subset B$ a non-empty subset
consisting of simplexes of dimension $\leq m$.  Then $|A|$ has the
homotopy type of a finite CW complex of
dimension
$\leq m$.

\end{lemma}

\noindent{\bf Proof:}  Let $\sigma$ be a top-dimensional simplex in $A$,
and let $C$ denote $\{\tau\in A| \tau\leq\sigma\}$.  Clearly $|C|$
contracts to the barycenter $\hat{\sigma}$ by sliding points along
``radial'' lines.  So, if $C=A$, we are done.  Otherwise, set
$A^{\prime}=A\setminus\{\sigma\}$ and
$C^{\prime}=C\setminus\{\sigma\}$. Note that
$|A^{\prime}|,\:|C|,\:\mbox{and}\:|C^{\prime}|$ are closed subsets of
$|A|$ satisfying $|A^{\prime}|\cup|C|=|A|$ and
$|A^{\prime}|\cap|C|=|C^{\prime}|$.  Moreover, by using the join
structure coming from barycentric subdivision, it is not hard to show
that $|C^{\prime}|$ is a neigborhood deformation retract (NDR) in both
$|A^{\prime}|$ and $|C|$\, so that the inclusion maps 
$|C^{\prime}|\hookrightarrow|A^{\prime}| $ and
$|C^{\prime}|\hookrightarrow|C|$ are cofibrations, (cf. [~], p. ). By
induction on the number of simplexes, both $|A^{\prime}|$ and
$|C^{\prime}|$ are homotopy-equivalent to finite CW complexes of
dimension $\leq m$ and $m-1$, respectively, and we have shown that this
also holds trivially for $|C|$. Standard homotopy-theoretic arguments now
yield the desired result for $|A|= |A^{\prime}|\cup|C|$. This completes
our proof.
\interpar

In order to apply the foregoing to the case of pairs, we use the
following result, whose proof is an easy exercise in homotopy theory and
will be omitted.
\interpar

\begin{lemma}

Suppose that a space $B$ has the homotopy type of an $\ell$-dimensional CW
complex, and a space $A$ has the homotopy type of a $k$-dimensional CW
complex.  Let $f:A\rightarrow B$ be a map with mapping cylinder $M(f)$.
Then the pair $(M(f),A)$ has the homotopy type of a CW pair $(L,K)$ such
that $dimL= max(k+1,\ell)$.

\end{lemma}

\noindent{\bf Completing the proof of
Proposition 21:}

\emph{Case a. $(T,S)=(P\setminus\Sigma,\emptyset)$.} In this case,
$P=|K|=|K^{k-1}|$, and
$\Sigma=|K^{k-3}|$. By Lemma 5.2, then, $P\setminus\Sigma$ has the
homotopy type of a graph, as desired.

\emph{Case b. $(T,S)=(Q\setminus\Sigma,P\setminus\Sigma)$.}  In this case,
$Q=|L|=|L^{k}|$ and $\Sigma\supseteq |L^{k-3}|$ (see the Remark at the
end of subsection 5.2).  Thus, Lemma 5.2 implies that $Q$ has the homotopy
type of a 2-dimensional CW complex. Since
$P\setminus\Sigma$ is the boundary of $Q\setminus\Sigma$, it posseses a
collar neighborhood. This implies that the pair $(Q,P)$ is homeomorphic
to the pair $(M(i),P)$, where $i$ is the inclusion $P\hookrightarrow Q$.
Now apply Case 1 (to $P$) and Lemma 5.4, to conclude that $(Q,P)$ has
the homotopy type of a 2-dimensional CW pair.

\emph{Case c. $(T,S)=(R\setminus\Sigma,Q\setminus\Sigma)$.}  In this
case $R=|N|=|N^{k+1}|$ and $\Sigma\supseteq |N^{k-3}|$.  Therefore,
$R\setminus\Sigma$ has the homotopy type of a 3-dimensional CW complex.
Arguing as in Case 2, we see that $(R,Q)$ has the homotopy type of a
3-dimensional CW pair.

This completes the proof of {Proposition 21. 
\interpar

\subsection{Smoothing singular and WPS circuits} 

In this section we show how to smooth singular circuits and bordisms
away from subcomplexes of codimension two. We continue with the
notation of the preceding subsections. 

\begin{prop}\label{prop:smoothmap}

Let $(Q,a)$ be a relative singular $k$-circuit in $(X,A)$, and suppose
that
$\alpha$ is a smoothing of $Q\setminus\Sigma$.  Suppose that
$a|(P\setminus\Sigma)_{\alpha}$ is smooth, and choose any
$\epsilon>0$.  Then there exists $(Q,b)$, a relative singular $k$-circuit
in
$(X,A)$, such that
\begin{enumerate}

  \item $b|(Q\setminus\Sigma)_{\alpha}$ is smooth.

  \item $b|\Sigma\cup P=a|\Sigma\cup P.$

  \item $d(a(x),b(x))<\epsilon$, for all $x\in Q$.

  \item $L(b|Q\setminus\Sigma)=b(\Sigma)$ and
$L(b|P\setminus\Sigma)=a(P\cap\Sigma)$.
\end{enumerate}

\end{prop}

\noindent{\bf Proof:}  Let $\delta:Q\rightarrow[0,\epsilon)$ be a
continuous function such that $\Sigma\cup P =\delta^{-1}(0)$.  Choose a
map
$b_{1}|(Q\setminus\Sigma)_{\alpha}\rightarrow X$ which is a smooth
extension of $a|(P\setminus\Sigma)_{\alpha}$ and a $\delta$-approximation
of
$a|(Q\setminus\Sigma)_{\alpha}$.  Define $b=(a|\Sigma)\cup
b_{1}$.  It is not hard to check that $b$ is continuous.  Properties (a)
- (c) are true by construction.  Assertion (d) is an immediate
application of Proposition 3i. This completes the proof.

\begin{corr}
Suppose that the relative singular $k$-circuit $(Q,a)$ in
Proposition 22 is WPS.  Then the resulting
$k$-circuit $(Q,b)$ is WPS, and
\begin{enumerate}
  \item $\ell d(b|(Q\setminus\Sigma)_{\alpha})\leq k-2.$
  \item $\ell d(b|(P\setminus\Sigma)_{\alpha})\leq k-3.$
\end{enumerate}
\end{corr}
{\bf Proof:}  That $(Q,b)$ is WPS when $(Q,a)$ is follows immediately from
the construction.  The assertions concerning limit dimension follow
immediately from assertion d of the proposition, together with the definition of the
singular set $\Sigma$ given above Proposition 20.  This completes the proof.
\interpar

The following variation will be useful later.

\begin{corr}
Let $(Q,a)$ be a relative singular $k$-circuit in $(X,A)$, and let
$\Sigma$ be as before.  Then there exists a smoothing $\alpha$ of
$Q\setminus\Sigma$ and a relative WPS $k$-circuit $(Q,b)$ in $(X,A)$
such that
\begin{enumerate}
  \item $b|(Q\setminus\Sigma)_{\alpha}$ is smooth.
	 \item $(Q,b)$ is bordant to $(Q,a)$
  \item $L(b|Q\setminus\Sigma)=b(\Sigma)$ and
$L(b|P\setminus\Sigma)=b(P\cap\Sigma)$.
\end{enumerate}
	
\end{corr}
{\bf Proof:}  First smooth $Q$ by applying
Proposition 21 a,b.  Denote the smoothing
by
$\alpha$. Set
$a\pri=a|P$ and apply Proposition 22 to $(P,a\pri)$ to get a
map $b\pri$ that $\epsilon$-approximates $a\pri$, and is
smooth on $(P\setminus\Sigma)_{\alpha}$. For $\epsilon$ small enough, we
can obtain
$b\pri$ homotopic to $a\pri$.  Consequently, by the homotopy extension
 property, $b\pri$ extends to a continuous map ${b\pri}\pri$ on $Q$, which
is homotopic to $a$ as a map of pairs. Now apply
Proposition 22 and Corollary 23 to $(Q,{b\pri}\pri)$ to
obtain a WPS map $b:(Q,P)\rightarrow(X,A)$ that extends $a\pri$,
approximates
${b\pri}\pri$, and is smooth on $(Q\setminus\Sigma)_{\alpha}$.  We choose
the approximation good enough to insure that $b$ is homotopic to
${b\pri}\pri$, hence homotopic to $a$.  The first assertion of the
corollary is now true by construction, and the third is true by the
argument previously given. Finally, the second assertion holds by Rourke
and Sanderson's theorem (Theorem 4), because $a$ and $b$ are homotopic,
implying  $a_{\ast}([Q])=b_{\ast}([Q])$.  This completes the proof.
\interpar

Closely similar arguments prove the next two results.  We
leave details to the reader.

\begin{prop}
Let a constant $\epsilon>0$ be given, and let $(R,c)$ be a singular
circuit null bordism in
$(X,A)$ of the singular relative
$k$-circuit $(Q,b)$.  Let $\Sigma$ be the exceptional or singular set for
$R$ as before.  Finally, assume that $\alpha$ is a smoothing of
$Q\setminus\Sigma$ such that $b|(Q\setminus\Sigma)_{\alpha}$ is smooth. 
Then there is a smoothing $\beta$ of $R\setminus\Sigma$ extending
$\alpha$ and a continuous map $d:R\rightarrow X$ such that
\begin{enumerate}

  \item $d|(R\setminus\Sigma)_{\beta}$ is smooth.
	 \item $d|\Sigma\cup Q=c|\Sigma\cup Q.$
  \item $d$ is an $\epsilon$-approximation of $c$.
  \item Any such $d$ necessarily satisfies
$$L(d|R\setminus\Sigma)=d(\Sigma)$$ and $$L(d|\overline{(\delta
R\setminus\Sigma)\setminus(Q\setminus\Sigma)})=c(\Sigma\cap\overline{\delta
R\setminus Q}),$$
where the first closure is relative to $\delta R\setminus\Sigma$ and the
second relative to $\delta R$.
\end{enumerate}
\end{prop}

\begin{corr}  Suppose that the nullbordism $(R,c)$ of the preceding
proposition is WPS.  Then the resulting nullbordism $(R,d)$ is WPS, and
\begin{enumerate}
  \item $\ell d(d|(R\setminus\Sigma)_{\beta})\leq k-1$
  \item $\ell
d(d|\overline{(\delta R\setminus\Sigma)_{\beta}\setminus
(Q\setminus\Sigma)_{\beta}})\leq k-2$.
\end{enumerate}
\end{corr}
\interpar

It is clear that the maps obtained in Corollary 24 and Corollary 26 are,
respectively a pseudocycle in $(X,A)$ and a nullbordism of pseudocycles
in $(X,A)$.  We are finally prepared, therefore, to construct the desired
homomorphism $\psi:H_{\ast}(X,A)\rightarrow\Ps{\ast}(X,A)$.

\section{ The natural equivalence $\psi:H_{\ast}\rightarrow\Ps{\ast}$}

Let $(X,A)$ be a smooth manifold pair in $\mathcal{C}^{prop}$, and choose
any class
$z\in H_{k}(X,A)$.  According to the theorem of Rourke and Sanderson
(Theorem 4, \S5), we may represent $z$ by a unique bordism class of
singular
$k$-circuits in $(X,A)$, say $(Q,a)$.  Choose any triangulation $(L,K)$
of $(Q,P)=(Q,\delta Q)$ and define the subpolyhedron $\Sigma$ of $Q$, all
as in the preceding subsections.  Corollary 24 produces a WPS $k$-circuit
$(Q,b)$ bordant to
$(Q,a)$, and a smoothing $\alpha$ of $Q\setminus\Sigma$ such that
$b|(Q\setminus\Sigma)_{\alpha}$ is smooth.  Further, Corollary 23 shows
that $((Q\setminus\Sigma)_{\alpha},b|)$ satisfies the limit dimension
constraints required of a pseudocycle. 

\begin{defn}
We define $\psi(z)$ to be the
\Ps{k}-bordism class of this pseudocycle.
\end{defn}

\begin{prop}
	$\psi$ is a well-defined homomorphism.  Furthermore, $\psi$ is
natural with respect to smooth proper maps of manifold pairs, and
it commutes with boundary maps of long exact sequences.
\end{prop}
\noindent{\bf Proof:} This follows from the results of \S5.  Thus, given
$z$ as above, different choices leading to the definition of
$\psi(z)$ would yield, say, a WPS relative $k$-circuit
$(\tilde{Q},\tilde{b})$ representing
$z$ and a corresponding smooth relative $\Ps{k}$-cycle
$((\tilde{Q}\setminus\tilde{\Sigma})_{\tilde{\alpha}},
\tilde{b}|(\tilde{Q}\setminus\tilde{\Sigma})_{\tilde{\alpha}})$. This
$k$-circuit is bordant to the original (as circuits).  Hence, we may use
Proposition 26 to obtain a smooth $\Ps{k}$-bordism between
$((Q\setminus\Sigma)_{\alpha},b|(Q\setminus\Sigma)_{\alpha})$ and
$((\tilde{Q}\setminus\tilde{\Sigma})_{\tilde{\alpha}},
\tilde{b}|(\tilde{Q}\setminus\tilde{\Sigma})_{\tilde{\alpha}})$, showing that $\psi(z)$ is
well-defined.  We leave details to the reader.  That $\psi$ is a
homomorphism follows immediately from the definitions. 

To see that $\psi$ is natural with respect to smooth maps of compact
manifold pairs, suppose that $f:(X,A)\rightarrow(Y,B)$ is such a map, and
let $z$ be as before.  Using the same notation as above, we have seen that
$z$ may be represented by a WPS relative $k$-circuit $(Q,b)$ such that the
restriction of $b$ to the smooth manifold $(Q\setminus\Sigma)_{\alpha}$ is
a smooth map of manifold pairs
$((Q\setminus\Sigma)_{\alpha},(P\setminus\Sigma)_{\alpha})\rightarrow(X,A)$.
This restriction represents
$\psi(z)$. By definition,
$((Q\setminus\Sigma)_{\alpha},f\circ(b|(Q\setminus\Sigma)_{\alpha}))$
represents
$f_{\ast}(\psi(z))$. Also by definition $(Q,f\circ b)$ represents
$f_{\ast}(z)$.  But $f\circ b|(Q\setminus\Sigma)_{\alpha}$ is a smooth map
of pairs, which, by the discussion above and the first part of this proof,
may be taken to represent
$\psi(f_{\ast}(z))$.  Therefore,
$f_{\ast}(\psi(z))=\psi(f_{\ast}(z))$, as required. 

The proof that $\psi$ commutes with boundary maps of long exact
sequences follows along similar lines.  Start with $z$ and $(Q,b)$ as in
the previous paragraph, obtaining $((Q\setminus\Sigma)_{\alpha},
b|(Q\setminus\Sigma)_{\alpha})$ representing $\psi(z)$.  Then
$(P,b|P)=(\delta Q,b|\delta Q)$ represents $\partial (z)$, and, by the first paragraph of
the proof,
$((P\setminus\Sigma)_{\alpha},b|(P\setminus\Sigma)_{\alpha})$ represents
$\psi(\partial(z))$.  On the other hand, this same pair clearly
represents
$\partial(\psi(z))$.  

This completes the proof of the proposition. 

\interpar

 The following corollary
summarizes some of what we have shown here and in Section 5.

\begin{corr}

$\psi:\Hm{\ast}\rightarrow\Ps{\ast}$ is a natural transformation of
functors on
$\mathcal{C}^{prop}$, and it induces maps of long exact sequences.

\end{corr}

\begin{prop}
$\psi:H_{\ast}(pt)\rightarrow\Ps{\ast}(pt)$ is the canonical
isomorphism described in Proposition 13.
\end{prop}

We leave this easy calculation to the reader.




\begin{prop}
$\psi$ is a natural equivalence of functors on $\mathcal{C}^{prop}_{hand}$.
\end{prop}
{\bf Proof:}  Begin by proving that $\psi$ is an isomorphism for the
following pairs: $(D^{k},\emptyset),\; (D^{k},S^{k-1}),\;
(S^{k},\emptyset),$ $(D^{k}\times D^{\ell},S^{k-1}\times D^{\ell})$. 
This parallels the usual inductive homology computation for these
pairs (cf. Proposition 14), using the Five Lemma at each stage. Since
both homology and pseudohomology map finite disjoint unions to direct
sums, $\psi$ is an isomorphism for all pairs that can be obtained from those listed by
taking finite disjoint unions.

Now choose
$(X,A)$ in
$\mathcal{C}^{prop}_{hand}$, where $X$ has dimension $n$.  By definition,
there is a nested sequence of codimension-$0$ submanifolds of $X$,
$\{X_{i}\},\; -1\leq i\leq n$,  such that, 
$ X_{-1}$ is a tubular neighborhood of $A$,  $X=X_{n}$, and each
$X_{i}$ smoothly deformation retracts to
$X_{i-1}$ together with finite number of disjointly attached $i$-handles. The
results of the previous paragraph, together with the excision and homotopy
properties for homology and pseudohomology, imply that $\psi$ is an
isomorphism for each pair
$(X_{i},X_{i-1})$. (Cf. the comment after Corollary 12 in \S4.)  The desired result then
follows inductively by  applying the Five Lemma to the long-exact homology and
pseudohomology sequences of the triples $(X_{i},X_{i-1},A)$.

\begin{prop}
Any natural transformation $H_{\ast}\rightarrow\Ps{\ast}$ on
$\mathcal{C}^{prop}$ that induces maps of long exact sequences and
equals the canonical isomorphism at $pt$ coincides with $\psi$ on
$\mathcal{C}^{prop}_{hand}$.
\end{prop}
{\bf Proof:} Let
$(X,A)$ and $\{X_{i}\}$ be as in the preceding proof. Define chain complexes ${\bf
C}^{\Psi}_{\ast}$ and ${\bf C}^{H}_{\ast}$ as follows: ${\bf
C}^{\Psi}_{i}(X,A)=\Psi_{i}(X_{i},X_{i-1})$; the boundary map $\partial_{i}$ is given by
the appropriate connecting homomorphism in the long exact $\Psi_{\ast}$-sequence for the triple
$(X_{i},X_{i-1},X_{i-2})$. Similarly for ${\bf C}^{H}_{\ast}$ . Let
${\bf H}^{\Psi}_{\ast}$ and ${\bf H}_{\ast}$ denote the corresponding
homology groups.  It is well known that ${\bf H}_{\ast}$ is naturally
equivalent to ordinary (singular) homology.  The same relation holds
between ${\bf H}^{\Psi}_{\ast}$ and \Ps{\ast}.  It will be convenient
for us to make this explicit, which we do in the case of $\Psi_{\ast}$.  A precise
analog holds for $\Hm{\ast}$.  So, consider the following sequence of
homomorphisms:

\begin{center}
\[\Ps{i}(X,A)\stackrel{j_{\ast}}{\rightarrow}\Ps{i}(X,X_{-1})\stackrel
{k_{\ast}}{\leftarrow}\Ps{i}(X_{i},X_{-1})\stackrel{\ell_{\ast}}
{\rightarrow}\Ps{i}(X_{i},X_{i-1})\stackrel{\partial}{\rightarrow}
\Ps{i-1}(X_{i-1},X_{i-2}),\] and
\[ker(\partial)\stackrel{\pi}{\rightarrow}{{\bf H}^{\Psi}}_{i}(X,A)\rightarrow
0,\]
\end{center}
where $j_{\ast}, k_{\ast}$ and $\ell_{\ast}$ are induced by inclusion, and
$\pi$ is the obvious projection. Note that $j_{\ast}$ is an isomorphism,
$k_{\ast}$ is onto, and $im(\ell_{\ast})\subseteq ker(\partial)$. In fact
$\partial=m_{\ast}\partial\pri$, where $m_{\ast}$ is the injection
$\Ps{i-1}(X_{i-1},X_{-1})\rightarrow\Ps{i-1}(X_{i-1},X_{i-2})$ induced by
inclusion, and $\partial\pri$ is the boundary homomorphism in the long exact
sequence of the triple $(X_{i},X_{i-1},X_{-1})$. It follows that
$im(\ell_{\ast})=ker(\partial)$.  Further,
$k_{\ast}$ factors through the inclusion-induced isomorphism
$\Ps{i}(X_{i+1},X_{-1})\rightarrow\Ps{i}(X,X_{-1}))$, so that
$$\ell_{\ast}(ker(k_{\ast})=\partial(\Ps{i+1}(X_{i+1},X_{i}))=B_{i}({\bf
C}^{\Psi}_{\ast}(X,A))\subseteq ker(\partial).$$

Together, these observations give a proof that the relation
$\pi\ell_{\ast}(k_{\ast})^{-1}j_{\ast}$ defines an isomorphism
\[\Ps{i}(X,A)\rightarrow {{\bf H}^{\Psi}}_{i}(X,A).\]

To deal with naturality, we introduce the notion of a \emph{handle-adapted} map.
Suppose that $(X,A)$ and $(Y,B)$ are smooth manifold pairs equipped with handle decompositions $\{X_{i}\}$ and $\{Y_{i}\}$, respectively.  If $i> dimX=m$, we set $X_{i}=X_{m}$, and similarly for $Y_{i}$. We then say that a smooth map of pairs $f:(X,A)\rightarrow(Y,B)$ is \emph{handle-adapted} (with respect to the given decompositions) if $f(X_{i})\subseteq Y_{i}$, for all $i$. A similar definition applies to homotopies. We make use of the following lemma, which is proved in virtually the same way as the analogue for cellular maps of CW complexes (e.g., see ~\cite{sp}, p. 404).

\begin{lemma} 
Suppose that $f$ is a smooth, proper map of smooth manifold pairs equipped with handle decompositions.  Then $f$ is properly and smoothly homotopic to a handle-adapted map.  Moreover, if two handle-adapted maps are properly and smoothly homotopic, then the homotopy may be chosen to be handle-adapted.
\end{lemma}

With the aid of this lemma, it is now easy to see that the chain complexes we have defined extend to functors on $\mathcal{C}^{prop}_{hand}$ and that the isomorphism $\Ps{i}(X,A)\rightarrow {{\bf H}^{\Psi}}_{i}(X,A)$ is natural with respect to smooth maps. Of course this is clear in the case of homology.

Now any natural transformation $H_{\ast}\rightarrow\Ps{\ast}$ that induces maps of long exact sequences as
hypothesized will induce a map of chain functors ${\bf C}^{H}_{\ast}\rightarrow{\bf
C}^{\Psi}_{\ast}$.  Moreover, if such a transformation equals the canonical
isomorphism at $pt$, it is easy to see that it must induce the same map of chain
functors as $\psi$.  Hence, it induces the same isomorphism of corresponding
homology functors ${\bf H}_{\ast}\rightarrow{\bf H}^{\Psi}_{\ast}$ as $\psi$.

The desired conclusion now follows from the fact that both the transformation in
question and $\psi$  commute with the homomorphisms
$j_{\ast},\;k_{\ast},\ell_{\ast},\;\partial$, and
$\pi$ used  above to define the isomorphisms
\[H_{\ast}(X,A)\rightarrow {\bf H}_{\ast}(X,A).\]
and
\[\Ps{\ast}(X,A)\rightarrow {\bf H}^{\Psi}_{\ast}(X,A).\]

This completes the proof.

\newpage

\appendix{\textbf{\LARGE Appendix on Corners}}
\vspace{.25in}

Corners appear naturally in
differential topology when, for example, one takes the product of
manifolds with non-empty boundaries, or, for another example, when a
manifold with non-empty boundary is separated by a codimension-one
submanifold that meets the boundary transversally. Most of the time, such
corners can be ``straightened'' or ``rounded'' on a case by case basis
without great inconvenience to the construction at hand. (See ~\cite{m3}, ~\cite{d}, or
the discussion below for a description of angle straightening.  Note, as described later,
we distinguish between straightening and rounding.) As a result, corners are rarely
considered in the literature (e.g., see ~\cite{k}, p.3).  However, there are contexts in
which a more systematic approach to corners is desirable, as for example when dealing with
spaces of embeddings.  J.Cerf gave the first such
account in ~\cite{c} for exactly this purpose. Some aspects of Cerf's treatment were
amplified and extended in the Cartan Seminar lectures of A. Douady ~\cite{d}.  Since that
time, the only other systematic treatment that I have been able to find is in the
unpublished e-manuscript of R. Melrose~\cite{me}.  

\section{Smooth manifold pairs}
We take as known all the usual
definitions involving smooth manifolds with boundary. For all integers
$m\geq 0$ and all finite sequences $I=(i_{1},\ldots,i_{k})$ such that
$0\leq k\leq m$ and $i_{1}<i_{2}<\ldots<i_{k}$---we allow the empty sequence
$\emptyset$, in which case $k=0$--- 
define
$R^{m}_{I}=\{x=(x_{1},\ldots,x_{m})\in R^{m}| x_{i}\geq 0$, for all
$i\in I\}$. We may use the notation $|I|$ to denote $k$. These are our model spaces. Note
that
$R^{m}=R^{m}_{\emptyset}$. If
$U$ is open in $R^{m}_{I}$ and
$h=(h^{1},\ldots,h^{n}):U\rightarrow R^{n}$ is a map, then we say that it
is smooth if each coordinate function $h^{i}$ extends to a $C^{\infty}$
function on an open subset of $R^{m}$.  This leads to the definitions of
diffeomorphisms, smooth charts, smooth atlases, and differentiable
structures in the usual way.  Note that if $V$ is open in $R^{n}_{J}$
and $h:U\rightarrow V$ is a diffeomorphism, then $m=n$, but 
$|I|$ is not necessarily equal to $|J|$ unless $0\in U$ and $h(0)=0$.  

\begin{defn}

\begin{enumerate}
  \item A smooth $m$-manifold with corners (or smooth $m$-manifold for
short) is a second-countable Hausdorff space $X$ together with a
differentiable structure $\{(U_{\alpha},h_{\alpha})\}$ for which the
image of each chart is an open subset of some $R^{m}_{I}$, $m$ fixed, $I$
variable.  Topologically, $X$ is an $m$-manifold with boundary.
  \item Let $X$ be a smooth $m$-manifold and $n$ a non-negative integer
$\leq m$.  A subset $A\subseteq X$ is called a smooth $n$-dimensional
submanifold of $X$ if every $a\in A$ is contained in a smooth chart
$(U,h)$ for $X$ such that $h(U\cap A)=h(U)\cap R^{n}_{J}$, for some
$J$.  Such an $A$ inherits a topology and a smooth atlas---hence a
differentiable structure---from $X$, making it into a smooth $n$-manifold.
  \item A smooth manifold pair $(X,A)$ is a topological pair such that
$A$ is a closed subset of $X$, $X$ and $A$ are smooth manifolds, and the
differentiable structure of $A$ is inherited from that of $X$. Note
that in this case the inclusion $A\hookrightarrow X$ is a proper, smooth
embedding.
\end{enumerate}
\end{defn}

\section{Corner indices and corner sets}

Fix the integer $m$, and choose any $p\in R^{m}_{I}$.  The
\emph{corner index} of
$p$ is the smallest integer $\ell$ such that $p$ has a neighborhood $U$ in
$R^{m}_{I}$ diffeomorphic to an open subset of $R^{m}_{J}$ with
$\ell=|J|$.  Thus, for example, the origin of $R^{m}_{I}$ has corner index
$|I|$ and $(1,1,\ldots,1)$ has corner index
$0$.  If $U$ and
$V$ are open subsets of $R^{m}_{I}$ and $R^{m}_{J}$, respectively, and
$h:U\rightarrow V$ is a diffeomorphism, it is easy to
show that
$h$ preserves the corner indices of points. It follows that we can
well-define the corner index of a point in a smooth manifold. We call the
maximal corner index of all points in $X$ the \emph{corner index} of
$X$. When there is no possibility of confusion, we shall omit the
adjective ``corner.'' It is easy to see that the corner index
of a point $p$ equals $|I|$ for any smooth coordinate chart $(U,h)$ around
$p$ such that $h(p)=0$ and $h(U)$ open in $R^{m}_{I}$.
\interpar

In general, the points of index $\geq
1$ in $X$ comprise the topological boundary of $X$, a closed set denoted
$bX$. The complement $X\setminus bX$ is called the \emph{manifold
interior of $X$} and is denoted $intX$. The points of index
$\geq 2$ comprise a closed subset
$cX\subseteq bX$,  called the corner set of $X$, (or the set of corners of
$X$).  $X\setminus cX$ and $bX\setminus cX$ are
smooth submanifolds of $X$ of index $\leq 1$, and clearly $b(X\setminus
cX)= bX\setminus cX$. The components of
$bX\setminus cX$ are called \emph{open faces} of $X$.  If $F$ is an open
face, it may happen that its closure $\overline{F}$ is
a smooth codimension-one submanifold of $X$, in which case we call
$\overline{F}$ a \emph{closed face} of $X$.  However, the example of one lobe of a
lemniscate ~\cite{c} (including its interior), shows that not every open face
is the interior of a closed face.   

When $X$
has index $0$, then $X$ is a standard smooth manifold with empty boundary.
When $X$ has index $1$, then it is a standard smooth manifold with
non-empty boundary. When
$X$ has index $2$, then $cX$ is a codimension-one submanifold of $bX$
with empty boundary and is a codimension-two smooth submanifold of $X$. 

Let $A$ be a connected, codimension-one submanifold of $X$ of corner index
$\leq 1$.  It is not hard to show that the intersection $intA\cap bX$ is
both closed and open in $intA$ and does not meet $cX$. Thus, either
$intA$ is contained in an open face of
$X$ or $intA\subseteq intX$. In the latter case, either $A\subseteq intX$
or $A$ meets $bX\setminus cX$ transversally in a union of components of
$bA$.   

\section{Straightening angles, introducing corners, and gluing}

The simplest way in which corners appear is in the product of two
manifolds with boundary $X\times Y$, in which case $b(X\times
Y)=bX\times Y\cup X\times bY$, and $c(X\times Y)= bX\times bY$. In this
case, the closure of every open face is a closed face.  For example,
when $X,\;Y,\;bX$, and $bY$ are connected, the closed faces of
$X\times Y$ are $bX\times Y$ and $X\times bY$.  
Using smooth collars
$bX\times [0,1)\rightarrow X$ and
$bY\times [0,1)\rightarrow Y$, we see that $c(X\times Y)$ has a
neighborhood $N$ diffeomorphic to $c(X\times Y)\times R^{2}_{2}$. We
identify $N$ with this product.  Milnor's construction of
\emph{straightening angles (or corners)}, ~\cite{m3}, ~\cite{d}, sends $N$
homeomorphically to
$c(X\times Y)\times R^{2}_{1}$ via the rule
$(x,\rho,\theta)\mapsto (x,\rho,2\theta)$, in which polar coordinates are
used for the second and third entries. This homeomorphism $h$ is then used
to endow
$N$ with a diferentiable structure, whereas
$X\times Y\setminus c(X\times Y)$ keeps the structure induced from that
of $X\times Y$.  These are compatible where they are jointly defined,
so they determine a differentiable structure, say $\alpha$, on $X\times
Y$. We denote $X\times Y$ with this structure $(X\times Y)_{\alpha}$ and
say that $(X\times Y)_{\alpha}$ is obtained from $X\times Y$ by
straightening angles (or corners).

Note that $h|\{p\in N| \theta(p)=0\}$ is a smooth embedding.  Similarly,
for $h|\{p\in N|\theta(p)=\pi/2\}$. Therefore, the closed faces of
$X\times Y$ inherit the same differentiable structures from $X\times Y$
as from $(X\times Y)_{\alpha}$.  

This procedure may be extended to all $X$ of index $\leq 2$ by making use
of the relevant tubular neighborhood theory (see ~\cite{d}).  The important point is that
the original inherited structures on
$X\setminus cX$, $bX\setminus cX$, and $cX$ remain unchanged, and
similarly for the structures on the closed faces of $X$.  All that
changes is the structure for $X$ at points on $cX$. 

The procedure of straightening corners can be reversed.  Thus,
suppose that $X$ is a smooth manifold with non-empty boundary and
(perhaps) corners. Choose a connected codimension-one, properly embedded
submanifold $B$
of $X$ contained in an open face of $bX$ and let $C$ be a component of $bB$.  Then $C$
has a tubular neighborhood in
$X$ with fibre
$R^{2}_{1}$.  As shown in ~\cite{d}, its differentiable structure can be
replaced by that of a corresponding tube with fibre $R^{2}_{2}$ by
``halving'' angles.  This \emph{introduces} a corner at $C$.  This
procedure is useful when gluing two manifolds together along a
codimension-0 submanifold of their boundaries, as we now describe in a
special case.

Suppose that $A$ and $B$ are two smooth $m$-manifolds and that $E$ and
$F$ are two diffeomorphic $(m-1)$-manifolds of index  $1$ which are
smooth, properly embedded submanifolds of $A$ and $B$, respectively, with
$E\subseteq bA$ and $F\subseteq bB$.  $E$
meets the corner set $cA$ in a disjoint union of components of $bE$, in which all
points have index $2$. 
Introduce corners at every component not meeting $cA$,calling the resulting manifold $A$
again.
 Do the same for $F$ and $B$. $bE$ and $bF$ become a union of components of
$cA$ and $cB$ with each point having index $2$ in $A$ and $B$, respectively. $E$
and $F$ become closed faces of $A$ and $B$, respectively.

Now in $A\sqcup B$, identify $E$ and $F$ via the given diffeomorphism,
calling the result $A+_{E=F}B$, or $A+B$, for short. Let $C$ denote the common image of
$E$ and $F$ in $A+B$, and let $D$ denote the common image of $bE$ and $bF$. It remains to
give
$A+B$ the structure of a smooth manifold compatible with the
structures of $A$ and
$B$.  For this, we make use of collar neighborhoods of $E$ and $F$ in $A$ and $B$ and 
quarter-disc tubular neighborhoods of
$bE$ and $bF$.  Restricting the collars to $intE$ and $intF$, these fit together in the
standard way to form a product neighborhood of $intC$ in $A+B$, which is given the product
smoothing.  The quarter-disc tubular neighborhoods fit together to form a half-disc
tubular neighborhood, which we use
to smooth a neighborhood of $D$.  Finally, use the structures on
$A\setminus E$ and $B\setminus F$ to smooth the complement
$A+B\setminus C$.  These all fit together to give a smooth
atlas for $A+B$.  The resulting smooth manifold contains $A$ and $B$
as smooth submanifolds.

\section{Rounding corners}

The process of straightening angles does not change the underlying
topological manifold, but when we are dealing with \emph{submanifolds}
with corners, it might be necessary to do so. We describe one method for
doing this.  Assume that $X$ is
compact for simplicity. Let
$A$ be a submanifold of $X$ of codimension zero and having corner index $\leq 2$. Let
$C$ be an open, connected subset of $cA$ over which the normal, orthogonal
$R^{2}_{2}$-bundle in $A$ is trivial.  Choose a fixed trivialization, and
use it to identify the normal bundle with $C\times R^{2}_{2}$.  Choose any
$\epsilon>0$, and let $\phi:[0,\infty)\rightarrow[0,\infty)$ be a smooth
function satisfying the following conditions:
\begin{itemize}
  \item $\phi\pri(t)<0$, for $t\in(0,\epsilon)$.
  \item $\phi^{(n)}(\epsilon)=0$, for all $n\geq 0$.
  \item $\phi(t)=0$,  for all $t\geq\epsilon$.
  \item $\phi(t)=\phi^{-1}(t)$, for $0\leq t\leq\epsilon$.
\end{itemize}

Now define $W_{\epsilon}=\{(x,s,t)\in C\times R^{2}_{2}$ such that
$t<\phi(s)\}$, regarded as a subset of $A$ via the previously mentioned
trivialization.
\emph{A priori} this set would appear to depend on the trivialization. 
However, as shown in ~\cite{d}, the group of the orthogonal normal bundle of
$cX$ is $Z_{2}$, generated by the automorphism of $R^{2}_{2}$ that
exchanges the factors.  This induces a self-bundle-map of $C\times
R^{2}_{2}$ that leaves $W_{\epsilon}$ invariant. Thus, given a component
$K$ of $cA$, we can cover its normal bundle by trivializations and then
piece together copies of the corresponding $W_{\epsilon}$'s that fit
together smoothly to obtain a subset of the  normal bundle that we call
$Y_{\epsilon}$. Set
$r(A)=A\setminus Y_{\epsilon}$.  This eliminates the corner of $A$ at $K$
and is called the `rounding' of
$A$ at
$K$.  Sometimes we delete reference to
$\epsilon$ or $K$, and sometimes we take $K$ to be a union of components
of $cA$. It is not hard to show, as in Lemma A.1, that there is a smooth
deformation retraction $(X,A)\rightarrow (X,r(A))$. Moreover $r(A)$ is
diffeomorphic to $A$ with angles straightened at $K$.

This procedure can be applied in the following situation. Let $f$ and $g$
be smooth, real-valued functions on $X$ such that: 
\begin{itemize}
  \item  $0$ is a regular value
of $f, g, f|bX,$\, and $g|bX$; 
  \item  $(0,0)$ is a regular value of $F=(f,g)$;
  \item  $F^{-1}(0,0)\cap bX=\emptyset$. 
\end{itemize}
Set $A=f^{-1}[0,\infty),
B=f^{-1}(-\infty,0], C=g^{-1}[0,\infty), D=g^{-1}(-\infty,0]$. Then both
$F^{-1}(R^{2}_{2})=A\cap C$ and
$F^{-1}(-R^{2}_{2})=B\cap D$ are smooth, codimension-$0$ submanifolds with
corners at $F^{-1}(0,0)$ in $intX$.  We can round these by the above
procedure. Furthermore, $\overline{X\setminus r(B\cap D)}$ is a
smooth submanifold of $X$ coinciding with $A\cup C$ outside a neighborhood
of
$F^{-1}(0,0)$ and deformation retracting onto $A\cup C$. By abuse of
notation we denote it by $r(A\cup C)$.

Now consider a handle
$H$ attached to a submanifold $W$ in $X$, supposing that $cW\subseteq bX$
and
$H\subseteq intX$. Let $g$ be a smooth real-valued function that
has regular value $0$ and $W= g^{-1}[0,\infty)$. Extend the handle $H$
slightly into $W$ and round the result, obtaining a smooth
codimension-$0$  manifold $K$ in $intX$ that equals $H$ outside of $W$ and
has boundary transverse to $bW$. Let
$f$ be a smooth real-valued function with regular value $0$ and
$K=f^{-1}[0,\infty)$. Then
$f$ and
$g$ have the properties above and
$W\cup K=W\cup H$. The manifold $r(W\cup K)$ is what we mean when we
speak of attaching the handle $H$ to $W$ (in $X$) and rounding the result.

If $A$ is a submanifold of $X$, and $B$ is a handle attached to $A$
in $X$, it is straightforward to show that $r(A\cup B)$ is diffeomorphic
to $A+B$. 

\section{Smoothing maps}

In this paper, corners are encountered in two settings: in domains of
smooth maps and in their targets.  We deal first with corners in the
domain.

\begin{lemma}
Suppose that we are given smooth manifolds $X$ and $Y$, a closed subset
$A$ of $X$, 
 a continuous map
$f:X\rightarrow Y$ that is smooth in a neighborhood of $A$, and a
continuous function
$\epsilon:X\rightarrow (0,1]$. Suppose further
that $Y$ is endowed with a metric.  Then there exists a smooth map
$g:X\rightarrow Y$ that is an
$\epsilon$-approximation of $f$ such that $g=f$ in a neighborhood of $A$
and $g=_{\infty}f$.
 \end{lemma}
{\bf Proof:}  Leaving aside the last condition asserted, the result is
completely standard in the case of smooth manifolds without corners.  The
standard proof carries over without essential change to the general
case.  It remains to show that the last condition can be satisfied. 
For this, we suppose that
$X$ is non-compact---since otherwise the assertion is vacuously
true---and we construct a continuous function $\delta:X\rightarrow (0,1]$
satisfying
$\delta=_{\infty}0$ and
$\delta\leq\epsilon$.  Then, we use the earlier statements in the lemma
to find a smooth $\delta$-approximation $g$ of $f$ 
that equals $f$ in a neighborhood of $A$. It is easy to check that
$g=_{\infty}f$. This completes the proof.
\interpar

It is useful to have a slight variant of this result.

\begin{cor}
Let $X$, $Y$, $A$, and $\epsilon$ be as in the lemma above, and
suppose additionally that $A$ is a codimension-zero submanifold of $X$. 
If $f:X\rightarrow Y$ is a continuous map such that $f|A$ is smooth,
then there exists a smooth map
$g:X\rightarrow Y$ that is an
$\epsilon$-approximation of $f$ such that $g|A=f|A$ and
$g=_{\infty}f$.
 \end{cor}
{\bf Proof:}  Results of ~\cite{c} imply that there is a smooth self-map $h$ of
$X$ that is arbitrarily close to the identity and retracts a neighborhood
of $A$ onto $A$.  The approximation may be chosen so close that $fh$ is
a $\delta/2$-approximation of $f$, where $\delta$ is as in the proof of
the above lemma.  Clearly
$fh$ is smooth in a neighborhood of $A$ and $fh|A=f|A$.  Apply the lemma
to find an
$\delta/2$-approximation $g$ of $fh$ that equals $fh$ in a neighborhood
of $A$.  It is easy to check that $g$ has the desired properties.  This
completes the proof.

\begin{cor}
Let $X$ and $Y$ be a smooth manifolds endowed with metrics, and suppose
that $X$ has corner index 2.  Let
$\alpha$ be a differentiable structure on $X$ obtained by smoothing
corners as described in the previous subsection.  Let $f:X\rightarrow Y$
be smooth and 
$\epsilon:X\rightarrow (0,1]$ any continuous function.  Then there is a
smooth
$g:X_{\alpha}\rightarrow Y$ such that
\begin{enumerate}
  \item $g$ is an $\epsilon$-approximation of $f$.
  \item $g=f$ outside any prescribed neighborhood of $cX$.
  \item $g=_{\infty}f$.
\end{enumerate}
\end{cor}
{\bf Proof:} Let $N$ be a prescribed neighborhood of $cX$, and let
$T\subseteq N$ be a closed tubular neighborhood of $cX$ in
$X_{\alpha}$.  Set $A=\overline{X\setminus T}$, and apply the previous
lemma. 
\interpar

\noindent{\bf Remark:} Note that item (c) above implies that $L(g)=L(f)$
(Lemma 2.1). 
\interpar

It takes a bit more work to prove the following relative version of this
corollary.

\begin{lemma}
We continue with the notation and hypotheses of Corollary E.3 and 
suppose additionally that $W$ is a disjoint union of closed faces of
$bX$.  Then the conclusion of Corollary E.3 holds, and furthermore, $g$
may be chosen to coincide with $f$ on $W$.
\end{lemma}
{\bf Proof (sketch):}  Since Corollary E.3 already gives the desired
result outside a neighborhood of $cX$ and the only possible  problem occurs near
$bW$, we may assume WLOG that $X$ is a neighborhood of $bW$, indeed, after corners are
smoothed, that $X=bW\times R^{2}_{1}$, with the notation chosen so that $W$ is identified
with
$bW\times[0,\infty)\times 0$. Let
$T$ be a tubular neighborhood of $cX$ as in the proof of Corollary E.3
above, and let $T^{\prime}$ be the union of those components of $T$ that
meet
$bW$, rounding corners so that the fibres of the bundle
projection $T^{\prime}\rightarrow bW$ are smooth discs which intersect
$bX$ in closed line segments forming a $1$-disc bundle over
$bW$. Let
$\delta:X\rightarrow(0,1)$ be an arbitrary, continuous function. Finally,
set $B=\overline{X\setminus T^{\prime}}\cup W$.

 We claim that there is  a smooth map
$h:X\rightarrow X$ with the following properties: 
\begin{itemize}
  \item $h=id_{X}$ on $B$.
  \item $h$ is a $\delta$-approximation of $id_{X}$. 
  \item $h^{-1}(B)$ is a neighborhood of $B$.
\end{itemize}

We leave the construction of such an $h$ to the reader. The function
$\delta$ may be chosen so small that
$f=_{\infty}fh$ and that $fh$ is an $\epsilon/2$-approximation of $f$. 
Moreover, $fh$ is smooth in a neighborhood of $B$ and coincides with $f$
on $B$.  We then apply Lemma E.1 to $fh$ to obtain a smooth
$\epsilon/2$-approximation
$g$ that equals $fh$ at infinity and coincides with $fh$---hence with
$f$---on a neighborhood of $B$.  $g$ is the desired
$\epsilon$-approximation of $f$.

This completes the sketch of the proof.
\interpar 

The following result is an immediate consequence of the foregoing.

\begin{cor}
Let $X$ be a smooth manifold of index $\leq 1$, let $A$ be a
codimension-zero submanifold of $X\setminus bX$, closed as a subset of
$X$, and let
$H:X\times[0,1]\rightarrow Y$ be a smooth map.  Smooth the corners of
$X\times[0,1]$ via a differentiable structure $\beta$, and endow
$X\times[0,1]$ and $Y$ with metrics.  Let
$\epsilon:X\times[0,1]\rightarrow(0,1]$ be a continuous function.
Then, there exists a smooth $\epsilon$-approximation of $H$,
$G:(X\times[0,1])_{\beta}\rightarrow Y$, which is equal to $H$ at
infinity and equals $H$ on $(X\times 0)\cup A\times[0,1]\cup (X\times 1)$.
\end{cor}

We now return to Corollary E.2, and we continue to use the notation of
that result and  its proof.  Let
$\beta$ be the differentiable structure on $X_{\alpha}\times [0,1]$
obtained by straightening the angles at the corner
$b(X_{\alpha})\times\{0,1\}$. As already remarked, the closed faces
$X_{\alpha}\times\{0,1\}$ and $bX_{\alpha}\times [0,1]$ keep their
original differentiable structures in
$(X_{\alpha}\times[0,1])_{\beta}$. 

We are interested in the relationship
between two close, smooth approximations of the map $f$.  It is a
standard fact in algebraic topology that when the approximation is close
enough, the two maps are homotopic, and this translates without difficulty
to differential topology.  Here, we take into account the smoothing of
corners, as well as what happens at infinity.

\begin{lemma}
Let $X,\;Y,\; A,\; f,\; \epsilon$ be as in Corollary E.2.  Then, we may
choose
$\delta\leq
\epsilon$ such that if $g_{0},g_{1}:X_{\alpha}\rightarrow Y$ are  smooth
$\delta$-approximations of
$f$ extending
$f|A$, there is a
smooth map $G:(X_{\alpha}\times[0,1])_{\beta}\rightarrow Y$ which
satisfies
\begin{enumerate}
  \item $G(x,i)=g_{i}(x)$, for all $x$ and for $i=0,1$.
  \item $G$ is a $2\delta\circ pr$-approximation of $f\circ pr$.
  \item $G|A\times[0,1]=f\circ pr|A\times[0,1]$. 
  \item $G=_{\infty}f\circ pr$.
\end{enumerate}
\end{lemma}
{\bf Proof (Sketch):} Corollary E.5 shows that we need only find a smooth
homotopy rel$A$, $X_{\alpha}\times[0,1]\rightarrow Y$, between $g_{0}$ and
$g_{1}$ that is suitably close to $f\circ pr$. There are two well-known
proofs of the classical analogue. One is inductive, proceeding stepwise
over a countable union of coordinate neighborhoods.  The other uses
differential geometry.  Both methods can be extended to the case of
manifolds with corners.  In particular, as is shown in ~\cite{d}, connections
can be defined for manifolds with corners in such a way as to respect
the corner structure.  Thus, geodesics are defined in the usual way, as
well as convex, normal geodesic neighborhoods.  Thus, two close
approximations to $f$ can be deformed to one another by sliding along
uniquely determined geodesic arcs.  Wherever the maps $g_{0}$ and
$g_{1}$ coincide, the homotopy is stationary.  This completes the sketch.

\interpar
\textsc{Department of Mathematics}

\noindent\textsc{Cornell University}

\noindent\textsc{Ithaca, NY}\, 14853

\noindent\texttt{\large kahn@math.cornell.edu}

\end{document}